\documentclass [11pt] {article}
\usepackage{epsfig}
\usepackage{amsmath,amssymb,amsthm}
\usepackage{enumerate}
\begin{document}
%
%
\newtheorem{theorem}{Theorem}
\newtheorem{prop}[theorem]{Proposition}
\newtheorem{lemma}[theorem]{Lemma}
\newcommand{\beqn}{\begin{equation}}
\newcommand{\eeqn}{\end{equation}}
\newcommand{\nn}{\nonumber}
%

\numberwithin{equation}{section}

\def\R{{\mathbb R}}
\def\L {{\cal L}}
\def\I{{\mathbb {\cal  I}}}
\def\cD{{\cal D}}
\def\K{{\cal K}}

\def \vb {{\mathbf v}}
\def \wb {{\mathbf w}}

\def \ub {{\mathbf u}}
\def \ubb {{\bar{\ub}}}

\def \fb {{\mathbf f}}
\def\fbb {{\bar{\fb}}}
\def\Fb{{\mathbf F}}

\def \gb {{\mathbf g}}
\def\gbb {{\bar{\gb}}}
\def\Gb{{\mathbf G}}

\def\Hb{{\mathbf H}}
\def\hb{{\mathbf h}}

\def\barf{{\bar{f}}}
\def\barg{{\bar{g}}}

\def\freg{{0 \leq \lambda x \leq t \leq \mu x}}
\def\greg{{0 \leq \mu x \leq t}}

\def\eb{{\mathbf e}}
\def\bp{{\mathbf p}}
\def\bq{{\mathbf q}}
\def\br{{\mathbf r}}
\def\bs{{\mathbf s}}

\def\abo{{\boldsymbol { \alpha}}}
\def\bbo{{\boldsymbol { \beta}}}
\def\pbo{{\boldsymbol {\phi}}}
\def\psibo{{\boldsymbol {\psi}}}
\def\ptbo{{ \boldsymbol{ \widetilde{\phi}}}}

\newcommand{\la}{\langle}
\newcommand{\ra}{\rangle}

\def \ut {{\widetilde{u}}}
\def \utb {{\mathbf \ut}}
\def \utbb {{\bar{\utb}}}

\def \ft {{\widetilde{f}}}
\def \ftb {{ \mathbf \ft}}
\def \ftbb {{ \bar{\ftb}}}

\def \gt {{\widetilde{g}}}
\def \gtb {{ \mathbf \gt}}
\def \gtbb {{ \bar{\gtb}}}

\def \at {{\widetilde{a}}}
\def \bt {{\widetilde{b}}}

\def\sone{{\sigma_1}}
\def\stwo{{\sigma_2}}

\def \tcoct {{ \frac{\tau}{\sone+\stwo}}}

\newcommand{\lt}{{\widetilde{l}}}

\newcommand{\cleq}{{~\preccurlyeq ~}}

\newcommand{\calL}{ {\cal L} }
\newcommand{\calT}{ {\cal T} }
\def\calB{{\cal B}}

\newcommand{\Ut}{ {\widetilde{U}}}
\newcommand{\At}{ {\widetilde{A}}}
\newcommand{\Bt}{ {\widetilde{B}}}

\def \Ah {{\hat{A}}}
\def \Bh {{\hat{B}}}
\def \vh {{\hat{\vb}}}

\newcommand{\tcx}{{t-cx}}
\newcommand{\tcox}{{t-\sone x}}
\newcommand{\tctx}{{t- \stwo x}}
\newcommand{\tlx}{ { t - \lambda x}}
\newcommand{\tmx}{ { t - \mu x}}

\def \dint {{\int_0^\infty \int_{-\infty}^\infty}}

\title{ Stability for an inverse problem for a two speed hyperbolic
pde in one space dimension}

\author{
Rakesh\\
Department of Mathematical Sciences\\
University of Delaware\\
Newark, DE 19716\\
~\\
Email: rakesh@math.udel.edu\\
\and
Paul Sacks\\
Department of Mathematics\\
Iowa State University\\
Ames, IA 50011\\
~\\
Email: psacks@iastate.edu }
\date{October 22, 2009}
\maketitle



\begin{abstract}
Suppose $A(x),B(x)$ are $2 \times 2$ matrices on an interval $[0,
\infty)$ and $C$ a constant diagonal matrix with distinct positive entries.
 Let $U(x,t)$ be the matrix solution of the
system of hyperbolic PDEs $ CU_{tt} - U_{xx} -A U_x - BU=0$ on
$[0,\infty)
\times \R$ with the initial condition $U(\cdot,t)=0$ for $t<0$ and the
boundary condition $U(0,t) = \delta(t)I_2$. We prove a stability
result for the inverse problem of recovering $A,B$ from
$U_x(0,\cdot)$. The solutions of the forward problem propagate with
two different speeds so techniques for inverse problems for a
single hyperbolic PDE are not applicable in any obvious way.
\end{abstract}


\section{Introduction}
Below, $\hspace{-0.1in} \cleq \hspace{-0.1in}$ will mean an
inequality up to a constant multiple, all functions will be real
valued, upper case letters such as $M$ will represent $2 \times 2$
matrices with entries $M_{ij}$, lower case bold letters such as
$\vb$ will represent $2 \times 1$ vectors with components $v_1,
v_2$. All convolutions will be in the $t$ variable if the
convolution involves a function of $x$ and $t$.
We define the operator $\calL$ by $
\calL \vb := C \vb_{tt} - \vb_{xx} - A\vb_x - B\vb
$ where $ C  = \left [
\begin{array}{rr}
\lambda^2 & 0 \\
0 & \mu^2
\end{array}
\right ]
$ with $0 < \lambda < \mu$, $A(x),B(x)$ are real valued $2 \times
2$ matrices and $\vb(x,t)$ is a $2 \times 1$ vector. Let $U(x,t)$
be the real valued matrix solution of the IBVP
\begin{align}
 \calL U = 0 \qquad &  \mbox{for} ~ (x,t) \in [0, \infty) \times \R
\label{eq:Upde} \\
U = 0 \qquad & \mbox{for} ~ t<0 \label{eq:Uic} \\
U (x =0,t) = \delta(t) I_2 \qquad & \text{for} ~t \in \R.
\label{eq:Ubc}
\end{align}
where $I_2$ is the $2 \times 2$ identity matrix.
We study the recovery of $A(\cdot), B(\cdot)$ or a subset of these
coefficients if we are given $U_x(x=0,\cdot)$.

Such an inverse problem arises in the examination of the
structural integrity of a composite beam; please see the
introduction of [MNS05] for a discussion of this application and
also for other references related to this application. The problem
of determining the (spatially varying) parameters for the
Timoshenko model of a beam (see [A73]), from measurements of the
deflection from the neutral axis and the twist in the
cross-section, also may be modeled as an inverse problem for a two
speed second order hyperbolic system on an interval with two
dependent variables; the entries of $A,B$ are made up of the
parameters in the Timoshenko beam model.

For future use we define $\ub,
\ubb$ to be the columns of $U$, that is $U=[\ub, \ubb]$.
Then $\ub$ and $\ubb$ also satisfy (\ref{eq:Upde}), (\ref{eq:Uic})
but satisfy the boundary condition
\beqn
\ub (x=0,t) = \delta(t) \eb_1,
\qquad
\ubb (x= 0,t) =  \delta(t) \eb_2
\label{eq:uubBC}
\eeqn
where $\eb_1$ and $\eb_2$ are the columns of $I_2$. There are two
speeds of propagation associated with $\calL$, namely $1/\lambda$
and $1/\mu$ and $u_1,
\bar{u}_1$ are the fast moving components and $u_2,
\bar{u}_2$ the slow moving components of $\ub$ and $\ubb$
respectively. It is this feature of the problem which makes it
difficult to apply any obvious modification of the inversion
schemes popular for inverse problems for a single hyperbolic PDE
in one space dimension.

In Theorem \ref{thm:progwave} we show that
(\ref{eq:Upde})-(\ref{eq:Ubc}) has a unique solution in $C^2([0,
\infty) , \cD'(\R))$ and we give a progressing wave expansion of
$U$. We postpone the statement of the theorem about the existence
and the structure of $U$ to the end of this section since the
statement is quite long and follows from the standard progressing
wave expansion technique; we want to draw attention to the more
interesting results in Theorems \ref{thm:stability} and
\ref{thm:slow} stated below.

Let $D=\text{diag}(A)$ be the diagonal matrix formed by taking
just the diagonal entries of $A$. Define  the diagonal matrix
$M(x) := e^{ - \frac{1}{2}
\int_0^x D(y) \, dy}$ and define $\vh := M^{-1} \vb$. Then we may show that
\begin{align*}
C \vb_{tt}  -  \vb_{xx} - A \vb_x - B \vb & = M( C \vh_{tt} -
\vh_{xx} - \Ah \vh_x - \Bh \vh)
\end{align*}
where $\Ah = M^{-1}(A-D)M$ and $\Bh= M^{-1}(D^2/4 - D'/2 - AD  +
B)M$. Note that the diagonal entries of $\Ah$ are zero.
Further, $M(0)=I$ so $U(0,t)=\hat{U}(0,t)=\delta(t)I$ and $
(M^{-1}U)_x(0,t)=  \hat{U}_x(0,t). $ So for every pair $(A,B)$ one
can construct a pair $(\Ah, \Bh)$ with the same data $(M^{-1}
U)_x(0,t)$ as $(A,B)$ except that the diagonal entries of $\Ah$
are zero. {\em Hence, below we will study only the situation where
the diagonal entries of $A$ are known.}

Define the operator $\calL^T$ by
\[
\calL^T \vb := C \vb_{tt} - \vb_{xx} + (A^T\vb)_x - B^T\vb
=C \vb_{tt} - \vb_{xx} + A^T\vb_x - (B-A')^T\vb.
\]
If $\vb(x,t)$ and $\wb(x,t)$ are $2\times 1$ vectors then one may
show that
\begin{align}
\vb^T  \calL \wb - (\calL^T \vb)^T \wb
= & (\vb^T C \wb_t - \vb_t^T C \wb)_t + (\vb_x^T \wb -  \vb^T
\wb_x -\vb^T A \wb)_x
 \label{eq:divform}
\end{align}
implying $\calL^T$ is the formal adjoint of $\calL$. Hence $\calL$
will be formally self-adjoint iff $A^T=-A$ and $B^T=B-A'$, that is
iff the diagonal entries of $A$ are zero and $B-B^T = A'$.

An analysis of the linearized inverse problem with the
linearization done around $A=0, B=0$ gives an indication of the
results one may expect for the inverse problem under
consideration. When $A=0, B=0$ the solution of
(\ref{eq:Upde})-(\ref{eq:Ubc}) is
\beqn
U(x,t) =
\begin{bmatrix}
\delta(\tlx) & 0 \\
0 & \delta(\tmx)
\end{bmatrix}.
\label{eq:Ulin}
\eeqn
Hence the linearized forward problem about the trivial background
is the solution of the IBVP
\begin{gather}
C (\delta U)_{tt} - (\delta U)_{xx} = (\delta A)U_x + (\delta B)U,
\qquad (x,t) \in [0,
\infty) \times \R
\label{eq:dUpde}\\
(\delta U)(0,t)=0, \qquad \delta U =0 ~ \text{for} ~ t<0.
\label{eq:dUic}
\end{gather}
Here we assume that $\text{diag}(\delta A)=0$.

Fix a $\tau>0$. We use (\ref{eq:divform}) with $\calL$
corresponding to $A=0, B=0$, $\vb(x,t) = U(x, \tau-t)$ and
$\wb(x,t)= (\delta U)(x,t)$.  Integrating this relation over the
region $[0,
\infty) \times \R$, integrating by parts and using (\ref{eq:dUic}), (\ref{eq:Ulin}),
we obtain
\begin{align}
 \int_0^\infty \int_{-\infty}^\infty U(x, \tau-t)^T ((\delta A)U_x + (\delta
B)U)(x,t) \, dt \, dx & =
\int_{-\infty}^\infty U(0, \tau-t)^T (\delta U)_x(0,t) \, dt.
\label{eq:lin}
\end{align}
Now, using (\ref{eq:Ulin}) in (\ref{eq:lin}) and integrating one
may show that
\begin{align*}
(\delta U_x)_{11} (0, \tau)  &= \frac{1}{2 \lambda} (\delta
B)_{11}(x_f(\tau)), ~~~ (\delta U_x)_{12}(0, \tau)
=\frac{1}{\lambda + \mu} (\delta B)_{12}(x_m(\tau)) -
\frac{\mu}{(\lambda + \mu)^2} (\delta A)_{12}'(x_m(\tau))
\\
(\delta U_x)_{21}(0, \tau) & = \frac{1}{\lambda + \mu} (\delta
B)_{21}(x_m(\tau)) - \frac{\lambda}{(\lambda + \mu)^2} (\delta
A)_{21}'(x_m(\tau)), ~~~ (\delta U_x)_{22}(0, \tau) = \frac{1}{2
\mu} (\delta B)_{22}(x_s(\tau))
\end{align*}
where
\[
x_f(\tau) = \frac{\tau}{2 \lambda}, ~~ x_m(\tau) =
\frac{\tau}{\lambda + \mu}, ~~ x_s(\tau) = \frac{\tau}{2 \mu}
\]
are the lengths probed from the origin, in time $\tau$, by a round
trip using two fast waves, a fast and a slow wave, and two slow
waves respectively.
\begin{figure}[h]
\begin{center}
\epsfig{file=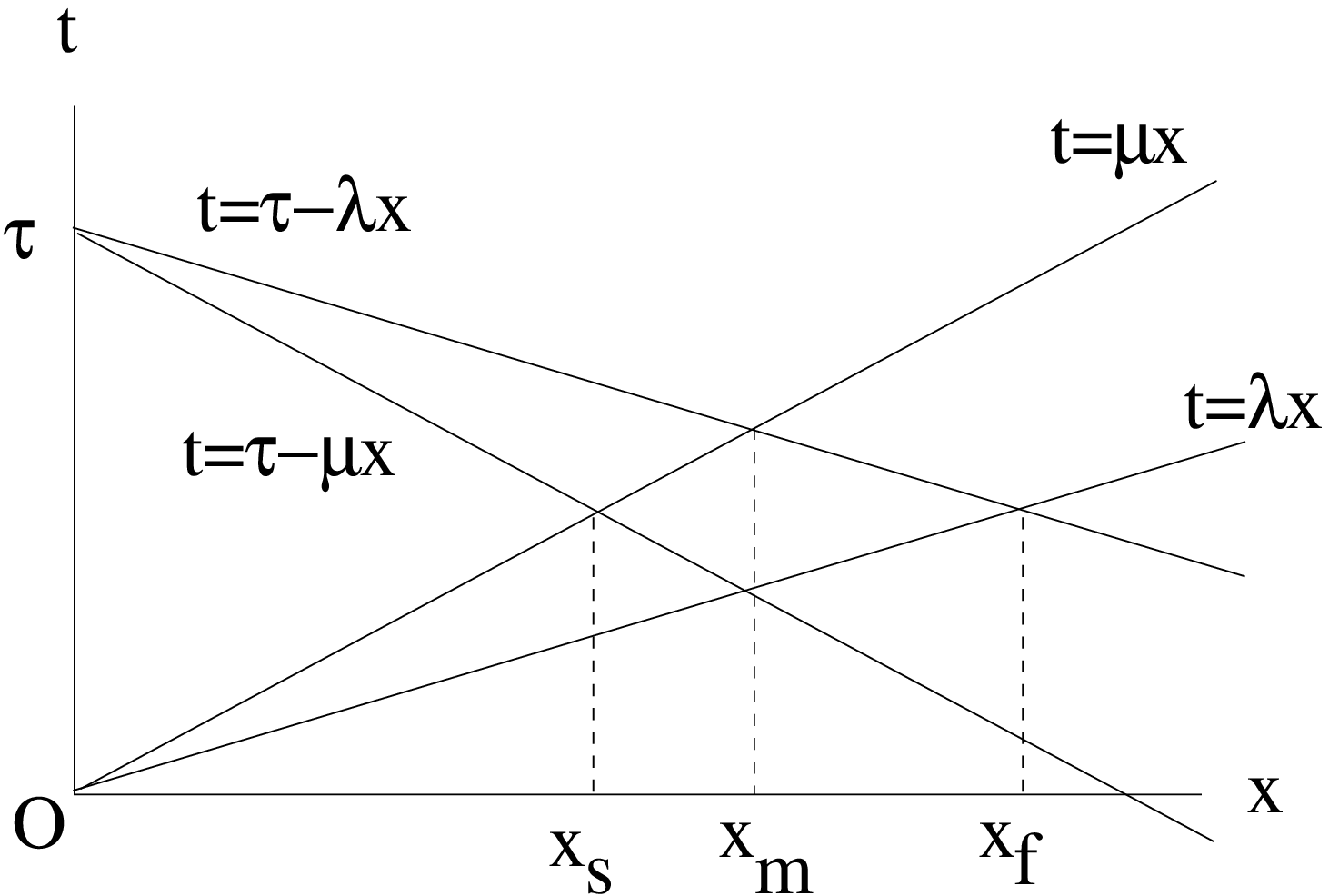, height=1.5in}
\end{center}
\caption{Definition of $x_s, x_m, x_f$}
\label{fig:xsxmxf}
\end{figure}

So, for this linearized inverse problem, where $(\delta U)_x(0,t)$
is given and $\delta A$, $\delta B$ are to be determined, one
recovers the combinations $\delta B_{11}$, $\delta B_{12} -
\frac{\mu}{\lambda + \mu}
\delta A_{12}'$, $\delta B_{21} - \frac{\lambda}{\lambda + \mu} \delta A_{21}'$ and
$ \delta B_{22}$. Hence if one is given two linearly independent
relations amongst $\delta B_{12}, \delta B_{21}, \delta A_{12}',
\delta A_{21}'$, which are independent of $\delta B_{12} -
\frac{\mu}{\lambda + \mu} \delta A_{12}'=0$ and $\delta B_{21} -
\frac{\lambda}{\lambda + \mu} \delta A_{21}'=0$,
 then one can recover $\delta A,
\delta B$ from $(\delta U)_x(0,\cdot)$. For example, if we are given the value
of $\delta A$ then one can recover $\delta B$. When the system is
self-adjoint we have $\delta B - (\delta B)^T =\delta A'$, that is
$\delta B_{12}- \delta B_{21} = \delta A_{12}' = - \delta
A_{21}'$. However, these relations are not independent of the two
relations mentioned above, so we need an additional relation or
the value of one of $\delta B_{12},
\delta B_{21}, \delta A_{12}', \delta A_{21}'$ would have to be part of the
data given.

This analysis suggests that, for the original inverse problem,
given $U_x(0,t)$ on an interval $[0, \tau]$ and the diagonal
entries of $A(\cdot)$, one may expect to recover only four out of
the remaining six coefficients in $A,B$, provided the other two
coefficients are given. Further, the values of these coefficients
will be recovered over intervals of different lengths which
suggest that there may be complications using the downward
continuation method popular for inverse problems for a single
hyperbolic PDE in one space dimension. However, if all the
coefficients except $b_{22}$ are known then there should be no
difficulty recovering $b_{22}$ with the use of a downward
continuation method.

Our main result is a stability result for the original inverse
problem and the proof reflects the discussion above. An
examination of the proof will show that one may prove stability in
more situations than covered in the statement of the theorem.
\begin{theorem}[Stability]\label{thm:stability}
Fix positive constants $X$ and $K$. Suppose $A, \At \in C^2[0,
X]$,
 $B, \Bt \in C^1[0, X]$ with
$\|A\|_{C^2} + \|B\|_{C^1} \leq K$, $\|\At\|_{C^2} + \|\Bt\|_{C^1}
\leq K$, and
\[
A(0)=\At(0), \qquad \text{diag}(A(\cdot)) =
0=\text{diag}(\At(\cdot)).
\]
Let $U$ and $\Ut$ be the solutions of
(\ref{eq:Upde})-(\ref{eq:Ubc})  corresponding to $A,B$ and $\At,
\Bt$ respectively, on the region $\{ (x,t) \, : \, 0 \leq x, ~ t +
\lambda x \leq \lambda X \}$. If either $A(\cdot)= \At(\cdot)$ or
the off-diagonal entries of $B(\cdot)$ and $\Bt(\cdot)$ are the
same, then
\begin{align}
|(B - \Bt)(x)| + |(A'-\At')(x)| \cleq \max_{[0, 2 \mu x]}
|U_x(0,\cdot) - \Ut_x(0,\cdot)|,
\qquad \forall x \in [0, \lambda X/\mu]
\label{eq:stability}
\end{align}
with the constant determined only by $X$, $K$, $\lambda$ and
$\mu$.
\end{theorem}
The theorem suggests that given $U_x(0,t)$ over the interval
$[0,T]$ one should be able to reconstruct (some of) the
coefficients over an interval $[0, T/(2 \mu)]$ - the interval
determined by the slower speed of transmission. Using the ideas
discussed earlier, one may derive a result similar to Theorem
\ref{thm:stability} if the hypothesis $\text{diag}(A) = 0 =
\text{diag}(\At)$ is replaced by the weaker hypothesis
$\text{diag}(A) = \text{diag}(\At)$.

For a $\bp(t) \in C^2(\R)$ with support in $[0, \infty)$, let
$\vb(x,t)$ be the solution of the IBVP
\begin{gather*}
\calL \vb = 0 ~ \text{in} ~ [0, \infty) \times \R\\
\vb(0,t) = \bp(t) ~ \text{for} ~ t \in \R,
\qquad \vb(\cdot,t)=0 ~\text{for} ~ t<0.
\end{gather*}
The fastest speed of propagation being $1/\lambda$, it is clear
that $\vb(x,t)$ will be supported in the region $0 \leq \lambda x
\leq t$. However, for certain choices of $\bp(\cdot)$, due to
cancelations, the support of $\vb(x,t)$ may lie in the slow region
$0 \leq \mu x \leq t$. In [BBI97] Belishev et al made an important
discovery where they showed that, if $\calL$ is formally
self-adjoint, then there is a unique function $l(t)$ (independent
of $\bp(t)$) so that if $p_1 = l*p_2$ then $\vb(x,t)$ is supported
in the slow region $0
\leq \mu x \leq t$. In fact, since $\vb = p_1*\ub + p_2*\ubb$
(with the convolution in $t$ alone), $l(t)$ is the unique function
so that $\ubb + l*\ub$ is supported in the slow region $0 \leq \mu
x \leq t$. Using some of the ideas in [BBI97], we have extended
their result to the general $A,B$ case and simplified the proof.
\begin{theorem}[Existence of slow waves]\label{thm:slow}
If $A \in C^2[0, \infty)$ and $B \in C^1[0, \infty)$ then there
exists a unique $l(\cdot)$ in $C^2[0,\infty)$ so that $\ubb(x,t) +
l(t)*\ub(x,t)$ is supported in the region $0 \leq \mu x \leq t$.
Further, for any $\tau>0$, $\|l\|_{C[0,
\tau(\mu-\lambda)/(\mu+\lambda)]}$ is bounded by a constant
determined only by $\lambda, \mu$ and $\|A\|_{C[0,x_f(\tau)]}$, $
\|B\|_{C[0,x_f(\tau)]}$.
\end{theorem}

In [BBI97] and [BI02], Belishev et al studied the inverse problem
considered in this article (for smooth coefficients though their
arguments are valid for less regular coefficients) except with the
additional requirement that $\calL$ be formally self-adjoint. In
this case there are only four coefficients to be determined but
then $U_x(0,\cdot )$ is also symmetric in this
case\footnote{For any $\tau>0$, using (\ref{eq:divform}) with
$\wb(x,t) = \ub(x,t)$ and $\vb(x,t) = \ubb(x,\tau-t)$ and
integrating over $[0,\infty) \times \R$ (using
$\wb(0,t)=\delta(t) \eb_1$, $\vb(0,t) =\delta(\tau-t)\eb_2$) we
obtain
\[
0 = \int_{-\infty}^\infty \eb_2^T \ub_x(0,t) \delta(\tau-t) -
\ubb_x^T(0,\tau-t) \eb_1 \delta(t) \, + \eb_2^T A \eb_1 \delta(t-\tau)
\delta(t)  \; dt = \partial_x U_{12}(0,\tau) -
\partial_xU_{21}(0,\tau).
\]
}
so the data $U_x(0,\cdot )$ consists of only three functions. With
this in mind, Belishev et al in [BBI97] and [BI02], for the
self-adjoint case, studied the recovery of $B$ from $U_x(0,\cdot)$
{\bf and} $l(\cdot )$. They showed that $B$ (and hence $A$) could
be reconstructed from $U_x(0,\cdot )$ and $l(\cdot )$. Further
(for the self-adjoint case) they characterized the range of the
map $(A(\cdot),B(\cdot)) \mapsto U_x(0, \cdot)$; they showed that
a function $r(t)$ is in the range of this map iff a certain
integral operator, defined in terms of $r(t)$, is positive
definite. Their proof showed that any pair of functions $(r(t),
l(t))$ defined over appropriate intervals, with  $r(t)$ satisfying
the ``positivity property'' is generated, in the above sense, by
some $A(\cdot), B(\cdot)$ associated with a self-adjoint $\calL$.

Since $l(\cdot )$ is not an experimentally measurable quantity,
in [BI03], again for the self-adjoint case, and assuming $A$ was
known, Belishev et al studied the recovery of $B$ (three unknown
quantities) from $U_x(0,\cdot )$. They showed that they could
reconstruct $l(t)$, at least over a small interval, and hence from
[BBI97] they could recover $B$ over a small interval. Using this
result Morassi et al in [MNS05] showed that if $A=0$ and $B$ is
symmetric (part of self-adjoint case) then the map $B
\mapsto U_x(0,\cdot )$ is injective (uniqueness in the inverse problem).
Our Theorem 1 covers the uniqueness (but not the reconstruction)
results in the above references and we provide a fairly simple
proof of stability for a more general situation. Belishev et al
use the Boundary Control Method which has proved effective for
reconstructions for several inverse problems for hyperbolic PDEs
and Morassi et al combine this with a downward continuation
argument in the frequency domain. We do not have a reconstruction
method even if $l(\cdot )$ is part of the data. Finally, [Ni91] is
a good starting point to read about the results of Nizhnik and his
school on inverse problems for two velocity systems.

Our proof of Theorem \ref{thm:stability} uses a trick similar to
the one used to analyze the linearized inverse problem above. This
trick was first used (as far as we know) in [SnSy88] for a single
hyperbolic PDE and then applied to a system of hyperbolic PDEs in
[Sa86], [SaSy87].

The existence and uniqueness of a weak solution of
(\ref{eq:Upde})-(\ref{eq:Ubc}) may be proved by appealing to
standard results but proving higher order piece-wise regularity
requires dealing with some quirks in two speed problems. The
following proposition characterizes the principal singularities in
$\ub$ and $\ubb$ and the existence theory associated with this
expansion.

\begin{figure}[h]
\begin{center}
\epsfig{file=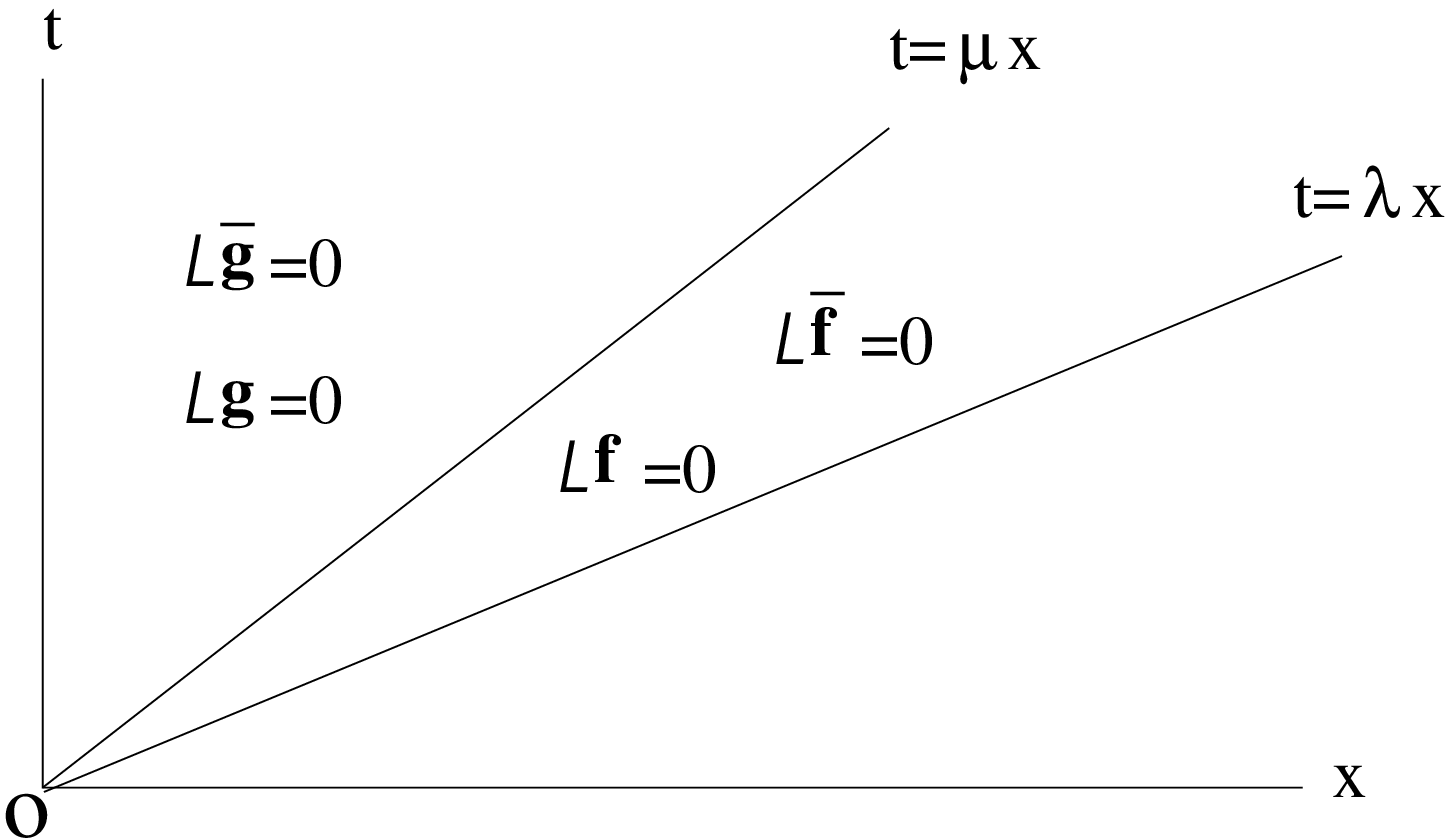, height=1.3in}
\end{center}
\caption{Domains of $\fb,\fbb,\gb,\gbb$}
\label{fig:fgfbgb}
\end{figure}

\begin{theorem}[Well posedness of the forward problem]\label{thm:progwave}
If $A \in C^2[0, \infty)$, $\text{diag}(A(\cdot ))=0$ and $B \in
C^1[0, \infty)$, then there exist unique solutions $\ub(x,t)$,
$\ubb(x,t)$ in $C^2([0,\infty), \cD'(\R))$ of (\ref{eq:Upde}),
(\ref{eq:Uic}), (\ref{eq:uubBC}). Further, for all $(x,t) \in [0,
\infty) \times \R$,
\begin{align}
\ub(x,t) & =  \delta(\tlx) \eb_1 + \fb(x,t) ( H(\tlx) - H(\tmx)) + \gb(x,t) H(\tmx),
\label{eq:ub}\\
\ubb(x,t) & =  \delta(\tmx) \eb_2 + \fbb(x,t) ( H(\tlx) - H(\tmx)) + \gbb(x,t) H(\tmx),
\label{eq:ubb}
\end{align}
where $\fb, \gb$ are $2 \times 1$ column vectors which are $C^2$
solutions of the characteristic IBVP (see Figure \ref{fig:fgfbgb})
\begin{align}
\calL \fb =0 ~~\text{in}~0 \leq \lambda x \leq t \leq \mu x,  \qquad \calL \gb =0 ~~ \text{in} ~
0 \leq \mu x \leq t,
\label{eq:fgpde}
\end{align}
with the boundary, characteristic and transmission conditions
\beqn
\gb(0,t) = 0, \qquad t \geq 0,
\label{eq:gbc}
\eeqn
\begin{align}
f_1(x, \lambda x) & =  \frac{1}{2 \lambda} \int_0^x b_{11}(z) \,
dz +
\frac{\lambda}{2(\mu^2 - \lambda^2)} \int_0^x a_{12}(z) a_{21}(z) \, dz,
\label{eq:f1} \\
f_2(x, \lambda x) &= \frac{ \lambda}{\lambda^2 -\mu^2} a_{21}(x),
~~
((\lambda^2 +\mu^2) f_{2t}  + 2 \lambda f_{2x} + \lambda a_{21}
f_1)(x,\lambda x) = b_{21}(x),
\label{eq:f2}
\end{align}
\begin{align}
(g_1 - f_1)(x, \mu x) &=0, ~~ (g_1 - f_1)_t(x, \mu x) = \frac{
\lambda \mu}{ (\mu^2 - \lambda^2)^2} a_{21}(0) a_{12}(x),
\label{eq:g1}\\
(g_2-f_2)(x, \mu x) &=  \frac{ \lambda}{\mu^2 -\lambda^2}
a_{21}(0).
\label{eq:g2}
\end{align}
Further $\fbb, \gbb$ are $C^2$ solutions of the characteristic
IBVP (see Figure \ref{fig:fgfbgb})
\begin{align}
\calL \fbb =0 ~~\text{in}~0 \leq \lambda x \leq t \leq \mu x  ,
\qquad \calL \gbb =0 ~~ \text{in} ~ 0 \leq \mu x \leq t,
\label{eq:fgbpde}
\end{align}
with the boundary, characteristic and transmission conditions
\beqn
\gbb(0,t) = 0, \qquad t \geq 0
\label{eq:gbbc}
\eeqn
\begin{align}
\barf_1(x, \lambda x) &= \frac{ \mu}{\lambda^2 - \mu^2} a_{12}(0),
\label{eq:barf1}
\\
\bar{f}_2(x,\lambda x) &=0,
~
\barf_{2t}(x, \lambda x) = \frac{\lambda \mu}{(\lambda^2 - \mu^2)^2} a_{12}(0)
a_{21}(x)
\label{eq:barf2}
\end{align}
\begin{align}
(\bar{g}_1-\bar{f}_1)(x, \mu x) &= \frac{\mu}{\mu^2 - \lambda^2}
a_{12}(x),
\label{eq:bg1f1} \\
(\barg_2 - \barf_2)(x, \mu x)  &= \frac{\mu}{2 (\lambda^2-\mu^2)}
\int_0^x a_{12}(z) a_{21}(z) \, dz + \frac{1}{2 \mu} \int_0^x b_{22}(z) \, dz,
\label{eq:bg2f2}
\\
 ( (\lambda^2+ \mu^2) (\barg_{1} - \barf_{1})_t & + 2 \mu(\barg_{1}
-
\barf_{1})_x  + \mu a_{12} (\barg_2 - \barf_2))(x, \mu x) = b_{12}(x).
\label{eq:bf1}
\end{align}
\end{theorem}
Using the ideas discussed earlier, one may derive a result similar
to Theorem \ref{thm:progwave} if the hypothesis $\text{diag} (A)
=0$ is dropped.

It would be reasonable to ask if results similar to Theorems 1,2,3
hold if the boundary condition (\ref{eq:Ubc}) is replaced by
$U_x(0,t) = \delta(t) I_2$ and for Theorem \ref{thm:stability} the
data is $U(0,t)$ instead of $U_x(0,t)$. We see no reason why the
same methods will not work after adjusting the order of the
singularity in $U(x,t)$, that is the most singular term in the
expansion of $U(x,t)$ would be $H(\tlx)$ and $H(\tmx)$ instead of
$\delta(\tlx)$ and $\delta(\tmx)$.

The rest of the paper consists of the following. In section
\ref{sec:stability} we prove Theorem \ref{thm:stability}.
In section
\ref{sec:slow} we prove Theorem \ref{thm:slow}.
Our proof uses some of the ideas in [BBI97] for the self-adjoint
case, but we do not use the Boundary Control Method machinery and
we think perhaps our proof is more transparent. In section
\ref{sec:progwave} we prove Theorem
\ref{thm:progwave} and Proposition \ref{prop:goursat} which is needed to complete
the proof of Theorem \ref{thm:progwave}. The proof of Theorem
\ref{thm:progwave} consists of two parts : a progressing wave
expansion and a well-posedness theory for a characteristic
transmission boundary value problem for a system of equations. The
progressing wave expansion part is standard but since the
expressions are not in the literature we give the expressions and
the derivation. The well-posedness theory for the characteristic
transmission boundary value problem for a system with two
velocities is not given in the literature though its proof uses
standard techniques except for the appearance of an unusual
transmission BVP problem for a single hyperbolic pde.

Finally we wish to thank Mikhail Belishev for discussions about
the problem considered in this article.


\section{Proof of Theorem \ref{thm:stability}}\label{sec:stability}
Extend  $A,\At$ as $C^2$ functions and $B, \Bt$ as $C^1$
functions, on $[0,\infty)$, with compact support, so that the
$C^2$ norms of $A, \At$ and the $C^1$ norms of $B, \Bt$, on $[0,
\infty)$, are bounded by a constant multiple of the corresponding
norms on $[0,X]$, with the constant independent of $A,\At, B,
\Bt$. Let $U=[\ub, \ubb]$ and $\Ut=[\utb, \utbb]$ be the solutions
of
 (\ref{eq:Upde})-(\ref{eq:Ubc}) corresponding to $A,B$ and $\At, \Bt$ respectively, over the
region $[0, \infty) \times \R$ guaranteed by Theorem
\ref{thm:progwave}. Further, let $l(\cdot )$ and $\lt(\cdot )$ be
the functions guaranteed by Theorem \ref{thm:slow} for the
operators corresponding to $A,B$ and $\At, \Bt$. Note that the
value of $U_x(0,t)$ and $\Ut_x(0,t)$ for $t \in [0, 2 \lambda X]$
is not affected by the extensions of $A, \At, B, \Bt$ because the
fastest speed of propagation is $1/\lambda$.

Define $\delta A := A - \At$, $\delta B := B - \Bt$, $\delta U :=
U - \Ut$, $(\delta a)_{ij} := a_{ij} - \at_{ij}$,  $(\delta
b)_{ij}  := b_{ij} - \bt_{ij}$,  and $|M| := \max_{ij} |m_{ij}|$.
Note that the diagonal entries of $\delta A$ are zero because of
the hypothesis. We will prove the stability by showing a Volterra
type estimate
\begin{align*}
| (\delta B - &  \frac{1}{\lambda + \mu}  (\delta A)'
\sqrt{C})(x)|
\\
& \cleq
\max_{[0, 2 \mu x]} |(\delta U)_x(0, \cdot )|  + \int_0^{2 \mu x} | (\delta U)_x(0,t)| \; dt
+ \int_0^x  | (\delta A)'(y)| + | ( \delta B)(y)| \, dy, ~~ x \in
[0, \lambda X/\mu]
\end{align*}
with the constant determined only by $\lambda, \mu$, $X$ and $K$.
Then Theorem \ref{thm:stability} follows from Gronwall's
inequality and the hypothesis that either $\delta A=0$ or the
off-diagonal entries of $\delta B$ are zero and the diagonal
entries of $\delta A$ are zero.

The progressing wave expansions of $\ub, \ubb$ are given by
(\ref{eq:ub}), (\ref{eq:ubb}) and from Theorem \ref{thm:progwave}
\begin{align}
\utb(x,t) &=\delta(\tlx) \eb_1 + \ftb (H(\tlx) - H(\tmx)) + \gtb
H(\tmx) \label{eq:utexp}\\
\utbb(x,t) &=\delta(\tmx) \eb_2 + \ftbb (H(\tlx) - H(\tmx)) + \gtbb
H(\tmx) \label{eq:utbexp}
\end{align}
with $\ftb$ and $\gtb$ having properties similar to $\fb, \gb$.
From Theorem \ref{thm:slow}, Theorem \ref{thm:progwave} and
Proposition \ref{prop:goursat} we have that the $C^0$ norm of $l,
\lt$ on any finite interval and the $C^2$ norms of $\fb, \ftb,
\gtb, \gtb$ on appropriate finite regions will be bounded by
functions of $\lambda, \mu, X$ and $K$ and parameters determining
the interval or the region. Since the regions of interest below
will be determined by $\lambda, \mu$ and $X$, one is assured that
all these norms are bounded by functions of $\lambda, \mu, X$ and
$K$.

We will use the following four pairs of vector functions
$\abo(x,t)$, $\bbo(x,t)$ defined on $[0,\infty) \times \R$ -
\begin{enumerate} [I.]
\item $\abo(x,t) = \ub(x,t) - \utb(x,t)$ and $\bbo(x,t) = \ub(x, \tau-t)$;
\item $\abo(x,t) = (\ubb - \utbb)(x,t) + \lt(t)*(\ub-\utb)(x,t)$ and
$\bbo(x,t) = \ubb(x,\tau - t) + (l*\ub)(x,\tau-t)$;
\item $\abo(x,t) = \ub(x,t) - \utb(x,t)$ and $\bbo(x,t) = \ubb(x,\tau - t) + (l*\ub)(x,\tau-t)$;
\item $\abo(x,t) = (\ubb - \utbb)(x,t) + \lt(t)*(\ub-\utb)(x,t)$ and $\bbo(x,t) = \ub(x, \tau-t)$.
\end{enumerate}
For each of these pairs we note that
\begin{itemize}
\item $\abo(0,t)=0$ on $\R$ and $\abo(x,t) =0$ for $t<0$
\item $\calL \bbo =0$ on $[0, \infty) \times \R$ and $\bbo(\cdot,t)=0$ for $t>>0$.
\end{itemize}
Hence using (\ref{eq:divform}) we have
\beqn
\int_0^\infty \int_{-\infty}^\infty \bbo^T \calL \abo \, dt \, dx
=
\int_0^\infty \int_{-\infty}^\infty \bbo^T \calL \abo - (\calL \bbo)^T \abo \, dt \, dx
= \int_{-\infty}^\infty \bbo(0,t)^T \abo_x(0,t) \, dt.
\label{eq:abiden}
\eeqn
In each of the four cases
\begin{align*}
\abo_x(0,t) &= (\ub - \utb)_x(0,t) ~ \text{or} ~ (\ubb - \utbb)_x(0,t) + \lt*(\ub - \utb)_x(0,t), \\
\bbo(0,t) & = \delta(\tau-t) \eb_1 ~ \text{or} ~ \delta(\tau-t) \eb_2 + l(\tau-t) \eb_1.
\end{align*}
Hence
\beqn
|\text{ RHS of (\ref{eq:abiden}) }|  \cleq |(U-\Ut)_x(0,\tau)| +
\int_0^\tau |(U-\Ut)_x(0,t)| \, dt
\label{eq:rhsest}
\eeqn
with the constant determined only by $X$ and $K$.

Estimating the LHS of (\ref{eq:abiden}), in each of the cases, may
involve one of the following estimates for $2 \times 1$ vectors
$\vb(x,t), \wb(x,t)$ which are $C^1$ and a continuous
 $2 \times2$ matrix $M(x)$. The derivation of these estimates is fairly straightforward with an integration by parts required for the first estimate.
\begin{align}
 \dint \vb^TM \, \wb  \delta'(\tau - \tcox) \, H(\tctx) \, dt \, dx  & \cleq
|M( \tcoct)| + \int_0^\tcoct |M(x)| \, dx;
\label{eq:dph}
\\
\dint \vb^TM \wb \, \delta(\tau - \tcox) \, \delta(\tctx) \, dt \, dx  & \cleq
|M( \tcoct)|;
\label{eq:dd}
\\
\dint \vb^TM \wb \, \delta(\tau - \tcox) \, H(\tctx) \, dt \, dx  & \cleq
 \int_0^\tcoct |M(x)| \, dx;
\label{eq:dh}
\\
\dint \vb^TM \wb \, H(\tau - \tcox) \, H(\tctx) \, dt \, dx  & \cleq
 \int_0^\tcoct |M(x)| \, dx,
\label{eq:hh}
\end{align}
with the constant determined only by the upper bounds on $|\vb|$,
$|\wb|$, $|\vb_t|$, $\wb_t|$  on the region $\{(x,t) \, :
\, 0 \leq
\stwo x
\leq t \leq \tau - \sone x \}$.

For future use we note that since $ \calL \delta U = (\delta A)
\Ut_x + (\delta B)
\Ut$ we observe that
\begin{align}
\calL (\ub - \utb) &= (\delta A) \utb_x + (\delta B) \utb,
\label{eq:uut} \\
\calL (\ubb - \utbb) &= (\delta A) \utbb_x + (\delta B) \utbb
\label{eq:uutb}.
\end{align}
Also, from the construction of $l(\cdot )$ and $\lt(\cdot )$ we
know that there are $C^1$ vector functions $\pbo(x,t)$ and
$\ptbo(x,t)$ so that
\begin{align}
(\ubb + l*\ub)(x,t)  &= \delta( \tmx) \eb_2 + \pbo(x,t) H(\tmx),
\label{eq:phi}\\
(\utbb + \lt*\utb)(x,t)  &= \delta( \tmx) \eb_2 + \ptbo(x,t)
H(\tmx)
\label{eq:phit}
\end{align}
and the $C^1$ norms of $\pbo$ and $\ptbo$ on appropriate finite
regions are bounded by $\lambda, \mu, X$ and $K$. For future use
we note that since $\utb$ is given by we may conclude that
\begin{align}
\utb_x(x,t) &= -\lambda \delta'(\tlx) e_1 + \ftb \delta( \tlx) + (\gtb - \ftb) \delta(\tmx) +
\ftb_x H(\tlx) + (\gtb_x - \ftb_x) H(\tmx).
\label{eq:utbx}
\end{align}
Below all constants are determined only by $\lambda, \mu, X$ and $K$.\\

\noindent{\bf Case I}\\
From (\ref{eq:ub}) we have
\[
\bbo(x,t) = \delta( \tau - \tlx) \eb_1 + \fb(x, \tau-t) H(\tau -
\tlx) + (\gb-\fb)(x, \tau-t) H(\tau-\tmx).
\]
From (\ref{eq:uut}) we note that $\calL \abo = (\delta A)
\utb_x + (\delta B) \utb$. Now $\utb, \utb_x$ are given by (\ref{eq:utexp}) and (\ref{eq:utbx}),
so some important contributions to the LHS of (\ref{eq:abiden})
from some singular terms in $\bbo^T \calL \abo$ are
\begin{align*}
\dint \eb_1^T (\delta A)\eb_1 \, \delta'(\tlx) \, \delta(\tau-\tlx)  \; dt \; dx &=0, \\
\dint \eb_1^T (\delta B) \eb_1 \delta( \tau-\tlx) \delta( \tlx) \; dt
\; dx &= (\delta B)_{11}(x_f(\tau))/(2 \lambda).
\end{align*}
All other terms on the LHS of (\ref{eq:abiden}) may be estimated
using (\ref{eq:dph})-(\ref{eq:hh}). Hence, using $(\delta
A)(0)=0$, we have
\begin{align}
|(\delta B)_{11}(x_f(\tau)) | &  \cleq | \text{RHS of
(\ref{eq:abiden})}| + |(\delta A)(x_f(\tau))| + \int_0^{x_f(\tau)}
|(\delta A)(x)| + |(\delta B)(x)| \, dx
\nn \\
& \cleq  |\text{RHS of (\ref{eq:abiden})}| +   \int_0^{x_f(\tau)}
|(\delta A)'(x)| + |(\delta B)(x)| \, dx.
\label{eq:pest}
\end{align}

\noindent{\bf Case II}\\
From (\ref{eq:uut}), (\ref{eq:uutb}) we note that $\calL \abo =
(\delta A)(
\utbb +
\lt*\utb)_x + (\delta B) (\utbb+ \lt*\utb)$
so
\begin{align}
\calL \abo (x,t) =  (\delta A) & \left ( -\mu \delta'(\tmx) \eb_2 - \mu  \ptbo(x,t) \delta( \tmx) + \ptbo_x(x,t) H(\tmx) \right ) \\
& + (\delta B) \left ( \delta( \tmx) \eb_2 + \ptbo(x,t) H(\tmx)
\right )
\label{eq:La}
\end{align}
and
\[
\bbo(x,t) = \delta( \tau - \tmx) \eb_2 + \pbo(x, \tau-t)
H(\tau - \tmx).
\]
Some important contributions to the LHS of (\ref{eq:abiden}) from
some singular terms in $\bbo^T \calL \abo$ are
\begin{align*}
- \mu \dint  \eb_2^T (\delta A) \eb_2 \delta'(\tmx) \delta(\tau-\tmx) &=0,\\
\dint \eb_2^T (\delta B) \eb_2 \delta(\tmx) \delta( \tau - \tmx) & =
\frac{1}{2 \mu} (\delta B)_{22}( x_s(\tau)).
\end{align*}
All other terms on the LHS of (\ref{eq:abiden}) may be estimated
using (\ref{eq:dph})-(\ref{eq:hh}). Hence, as before, we have
\begin{align}
|(\delta B)_{22}(x_s(\tau))| \cleq |\text{RHS of
(\ref{eq:abiden})}| +
\int_0^{ x_s(\tau)} |(\delta B)(x)| + |(\delta A)'(x)| \, dx.
\label{eq:sest}
\end{align}

\noindent
{\bf Case III}\\
From (\ref{eq:phi}) we have
\[
\bbo(x,t) =\delta( \tau - \tmx) \eb_2 + \pbo(x, \tau-t)
H(\tau - \tmx)
\]
and from (\ref{eq:uut}) we have $\calL \abo=(\delta A)\utb_x +
(\delta B)\utb$. Using $\utb, \utb_x$ given by (\ref{eq:utexp}),
(\ref{eq:utbx}), some important contributions to the LHS of
(\ref{eq:abiden}) from some singular terms in $\bbo^T \calL \abo$
are
\begin{align*}
-\lambda \dint \eb_2^T (\delta A)(x) \eb_1 \delta(\tau -\tmx)
\delta'(\tlx) &= -\frac{\lambda}{(\lambda + \mu)^2} (\delta A)_{21}'(x_m(\tau)),\\
 \dint \eb_2^T (\delta B) \eb_1 \delta(\tlx) \delta(\tau - \tmx) &=
 \frac{1}{\lambda + \mu} (\delta B)_{21}( x_m(\tau)).
\end{align*}
All other terms on the LHS of (\ref{eq:abiden}) may be estimated
using (\ref{eq:dph})-(\ref{eq:hh}). Hence, as before, we have
\begin{align}
|((\delta B)_{21}-\frac{\lambda}{\lambda + \mu} (\delta A)_{21}'
)(x_m(\tau)) |
\cleq |\text{RHS of
(\ref{eq:abiden})}| + \int_0^{ x_m(\tau)} |(\delta B)(x)| +
|(\delta A)'(x)| \, dx.
\label{eq:rest}
\end{align}

\noindent
{\bf Case IV}\\
From (\ref{eq:ub}) we have
\[
\bbo(x,t) =\delta( \tau - \tlx) \eb_1 + \fb(x, \tau-t)
H(\tau - \tlx) + (\gb - \fb)(x, \tau-t) H(\tau - \tmx)
\]
and $\calL \abo$ is the same as in Case II and is given by
(\ref{eq:La}). So some important contributions to the LHS of
(\ref{eq:abiden}) from some singular terms in $\bbo^T \calL \abo$
are
\begin{align*}
-\mu \dint \eb_1^T (\delta A)(x) \eb_2 \delta(\tau -\tlx)
\delta'(\tmx) &= \frac{-\mu}{(\lambda + \mu)^2} (\delta A)_{12}'(x_m(\tau)),\\
 \dint \eb_1^T (\delta B) \eb_2 \delta(\tmx) \delta(\tau - \tlx) &=
 \frac{1}{\lambda + \mu} (\delta B)_{12}( x_m(\tau)).
\end{align*}
All other terms on the LHS of (\ref{eq:abiden}) may be estimated
using (\ref{eq:dph})-(\ref{eq:hh}). Hence, as before, we have
\begin{align}
|((\delta B)_{12} - \frac{\mu}{\lambda + \mu} (\delta
A)_{12}'(x_m(\tau)) |
\cleq |\text{RHS of
(\ref{eq:abiden})}| + \int_0^{ x_m(\tau)} |(\delta B)(x)| +
|(\delta A)'(x)| \, dx.
\label{eq:qest}
\end{align}

Fix an $x$ in $[0,\lambda X/\mu]$ and define $t_s(x) = 2 \mu x, \,
t_m(x)=(\lambda + \mu)x, \, t_f(x)= 2 \lambda x$ to be the two-way
travel time to probe a distance $x$ at slow, mixed or fast speeds
respectively. Then (\ref{eq:pest}), (\ref{eq:sest}),
(\ref{eq:rest}), (\ref{eq:qest}), together with (\ref{eq:rhsest})
may be combined into
\begin{align}
|   ( (\lambda + \mu)(\delta B) -  (\delta A)' \sqrt{C} )(x)  |
 \cleq  & |(\delta U)_x(0,t_s(x))| + |(\delta U)_x(0,t_m(x))| +
|(\delta U)_x(0,t_f(x))|
\nn
\\
& + \int_0^{t_s(x)} |(\delta U)_x(0,t)| \; dt + \int_0^x |(\delta
A)'(z)| + |(\delta B)(z)| \; dz
\nn
\\
\cleq & \max_{[0, t_s(x)]} (\delta U)_x(0,\cdot ) + \int_0^x |(\delta
A)'(z)| + |(\delta B)(z)| \; dz.
\nn
\end{align}
\hfill{\bf QED}


\section{Proof of Theorem \ref{thm:slow}}\label{sec:slow}
Below all convolutions will be convolutions in the time variable
only. Because of the ideas discussed in the introduction it is
enough to prove Theorem
\ref{thm:slow} for the special case when $\text{diag}(A) = 0$ - we will assume that for the rest of the proof.

\begin{figure}[h]
\begin{center}
\epsfig{file=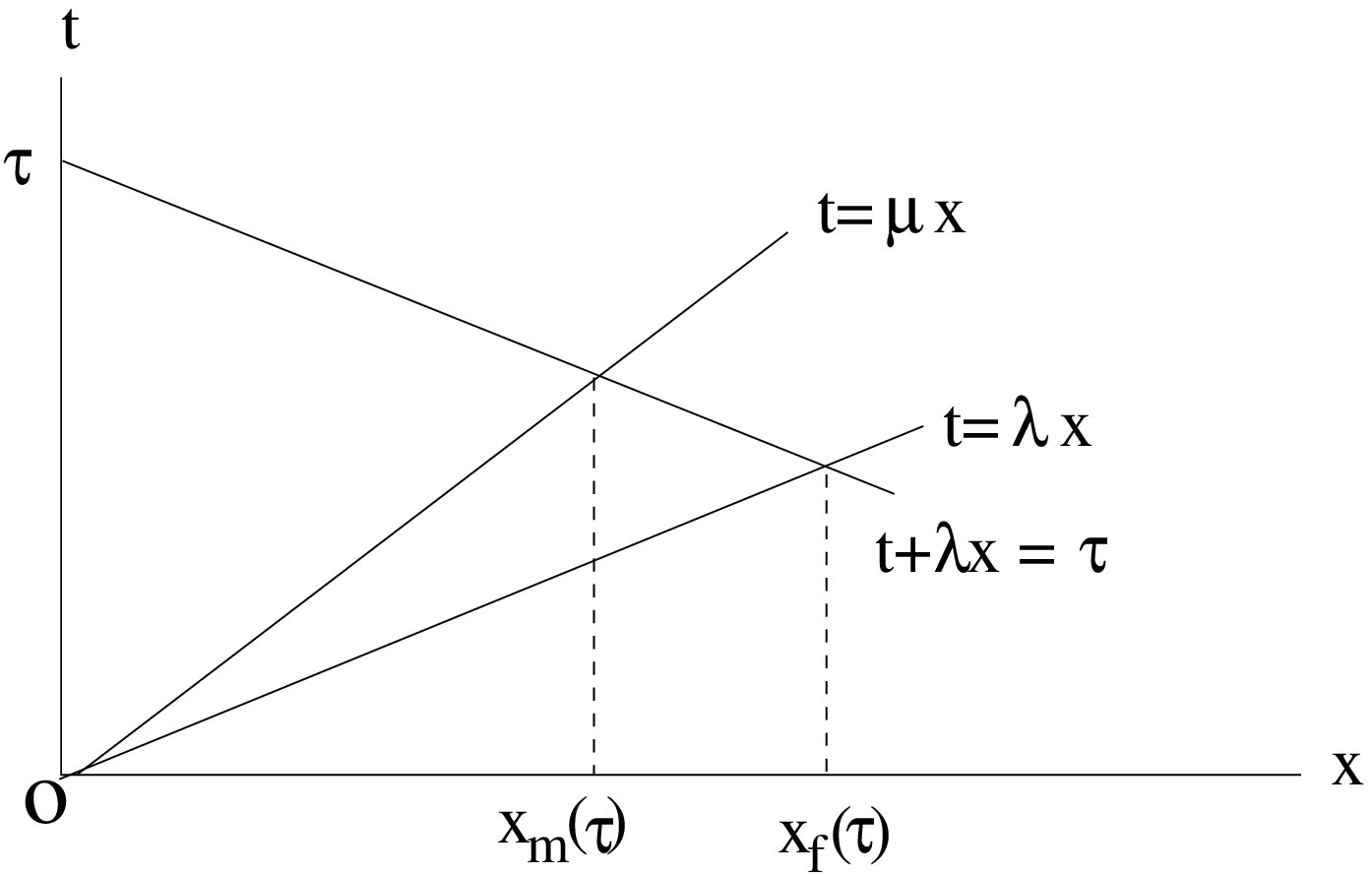, height=1.6in}
\end{center}
\caption{Interval of dependence for $l$}
\label{fig:Ypic}
\end{figure}

We must find an $l(t)$ supported in $[0,\infty)$ so that $
\vb(x,t) :=
\ubb(x,t) + l(t)*
\ub(x,t)$ is zero on $0<\lambda x \leq t < \mu x$. From
(\ref{eq:ub}), (\ref{eq:ubb}) we see that the most singular term
in $\vb(x,t)$ is $\delta(\tmx) \eb_2$ but this has no impact in
the region $0 < \lambda x \leq t <
\mu x$. So, for the rest of the proof we will identify $\vb(x,t)$
with $\vb(x,t) - \delta(\tmx) \eb_2$ over the region $0 \leq
\lambda x \leq t \leq \mu x$. Now over the region $0 < \lambda x
\leq t < \mu x$ one may observe that
\begin{align}
\vb(x,t) &= \fbb(x,t) + l(t)*\left ( \delta(\tlx) \eb_1 + \fb(x,t) H(\tlx) \right )
\nn \\
&= \fbb(x,t) + l(t-\lambda x) \eb_1 + \int_0^{t-\lambda x} l(s)
\fb(x, t-s) \, ds.
\label{eq:vbexp}
\end{align}
Hence we have to find an $l(t)$ so that
\begin{align}
l(t-\lambda x) + \int_0^{t-\lambda x} l(s) f_1(x, t-s) \; ds +
\barf_1(x,t)   & = 0, ~~0 \leq \lambda x \leq t \leq \mu x,
\label{eq:lone}
\\
\int_0^{t-\lambda x} l(s) f_2(x, t-s) \; ds
+ \barf_2(x,t)   & = 0, ~~0 \leq \lambda x \leq t \leq \mu x.
\label{eq:ltwo}
\end{align}
Fix a $\tau>0$; then rewriting (\ref{eq:lone}) for points on the
line $t+ \lambda x = \tau$, we seek a function $L(\cdot )$ so that
\beqn
L(\tau- 2\lambda x) + \int_0^{\tau- 2\lambda x} L(s) f_1(x,
\tau -\lambda x -s)
\, ds + \barf_1(x, \tau - \lambda x) =0, ~~ x_m(\tau) \leq x \leq x_f(\tau).
\label{eq:volterra}
\eeqn
The Volterra equation (\ref{eq:volterra}) has a unique solution
$L(\cdot )$ in $C[0, \tau(\mu-\lambda)/(\mu+ \lambda)]$. Since
$\fb$ and $\fbb$ are in $C^2$ and $L$ is continuous,
(\ref{eq:volterra}) implies that $L \in C^1[0,
\tau(\mu-\lambda)/(\mu+ \lambda)]$ which again by
(\ref{eq:volterra}) implies that $L \in C^2[0,
\tau(\mu-\lambda)/(\mu+ \lambda)]$. The $L(\cdot )$ constructed
depends on $\tau$ but we have to find an $L(\cdot )$ independent
of $\tau$ - except for the domain of $L$ which will depend on
$\tau$. Moreover this $L(\cdot )$ must also satisfy
(\ref{eq:ltwo}). Both these goals will be achieved if we can show
that $v_1(x,\tau - \lambda x)=0$ for $x_m(\tau) \leq x \leq
x_f(\tau)$ implies that $\vb(x,t)=0$ for $0 \leq \lambda x \leq
t\leq \mu x$, $t+ \lambda x \leq \tau$; see Figure \ref{fig:Ypic}.
Note that the supremum of $L(\cdot )$ on
$[0,\tau(\mu-\lambda/(\mu+\lambda))]$ is bounded above by a
function of the supremum of $f_1(x,t)$ and $\bar{f}_1(x,t)$ on the
region $0 \leq \lambda x
\leq t \leq\tau - \lambda x$. Hence, by Theorem \ref{thm:progwave}, the supremum of $L(\cdot )$
on $[0, \tau(\mu-\lambda)/(\mu+\lambda)]$ is bounded by a function
of the supremum
 of  $A(\cdot ),B(\cdot )$ on $[0, x_f(\tau)]$.

Below we use $L(t)$ instead of $l(t)$. From (\ref{eq:volterra})
and (\ref{eq:barf1}) we observe that
\[
L(0)= -\barf_{1}(\tau/\lambda, \tau) = \frac{-\mu}{\lambda^2 -
\mu^2} a_{12}(0)
\]
so using (\ref{eq:vbexp}) and (\ref{eq:barf1}), (\ref{eq:barf2}),
we see that
\[
\vb(x, \lambda x) = \fbb(x,\lambda x) + L(0) \eb_1 = 0.
\]
Further, from (\ref{eq:vbexp}), over $0 \leq \lambda x \leq t \leq
\mu x$ we have
\[
\vb_t(x,t) = \fbb_t(x,t) + L'(\tlx) \eb_1 + L(\tlx) \fb(x, \lambda x) + \int_0^{\tlx} L(s)
\fb_t(x,t-s) \, ds
\]
so using (\ref{eq:f2}), (\ref{eq:barf2}) we have
\begin{align*}
v_{2t}(x, \lambda x) & = \barf_{2t}(x, \lambda x) + L(0) f_2(x,
\lambda x) \\
& =
\frac{\lambda \mu}{ (\lambda^2-\mu^2)^2} a_{12}(0) a_{21}(x)
- \frac{\mu}{\lambda^2 - \mu^2} a_{12}(0)
\frac{\lambda}{\lambda^2- \mu^2} a_{21}(x) =0.
\end{align*}
Hence from (\ref{eq:ub}), (\ref{eq:ubb}) we see that
\begin{gather}
C \vb_{tt} - \vb_{xx} - A \vb_x - B \vb = 0, \qquad 0 \leq \lambda
x \leq t \leq \mu x, ~ t + \lambda x \leq \tau,
\label{eq:vlpde}\\
v_1(x, \lambda x) =0, ~~ v_2(x, \lambda x)=0, ~~ v_{2t}(x, \lambda
x)=0, ~~ 0 \leq x
\leq x_f(\tau).
\label{eq:v1char}
\end{gather}
We have to show that if $v_1(x,\tau - \lambda x)=0$ for $x_m(\tau)
\leq x \leq x_f(\tau)$, that is $v_1=0$ on $RS$ then $\vb(x,t)=0$
on the region ORS. This will follow from some energy identities -
the only complication being the two velocities. One could also do
this by setting $v_1=0$ on the relevant part of $t=\tau$ instead
of the $t+ \lambda x = \tau$ but one would not obtain the optimal
interval of dependence results.
\begin{figure}[h]
\begin{center}
\epsfig{file=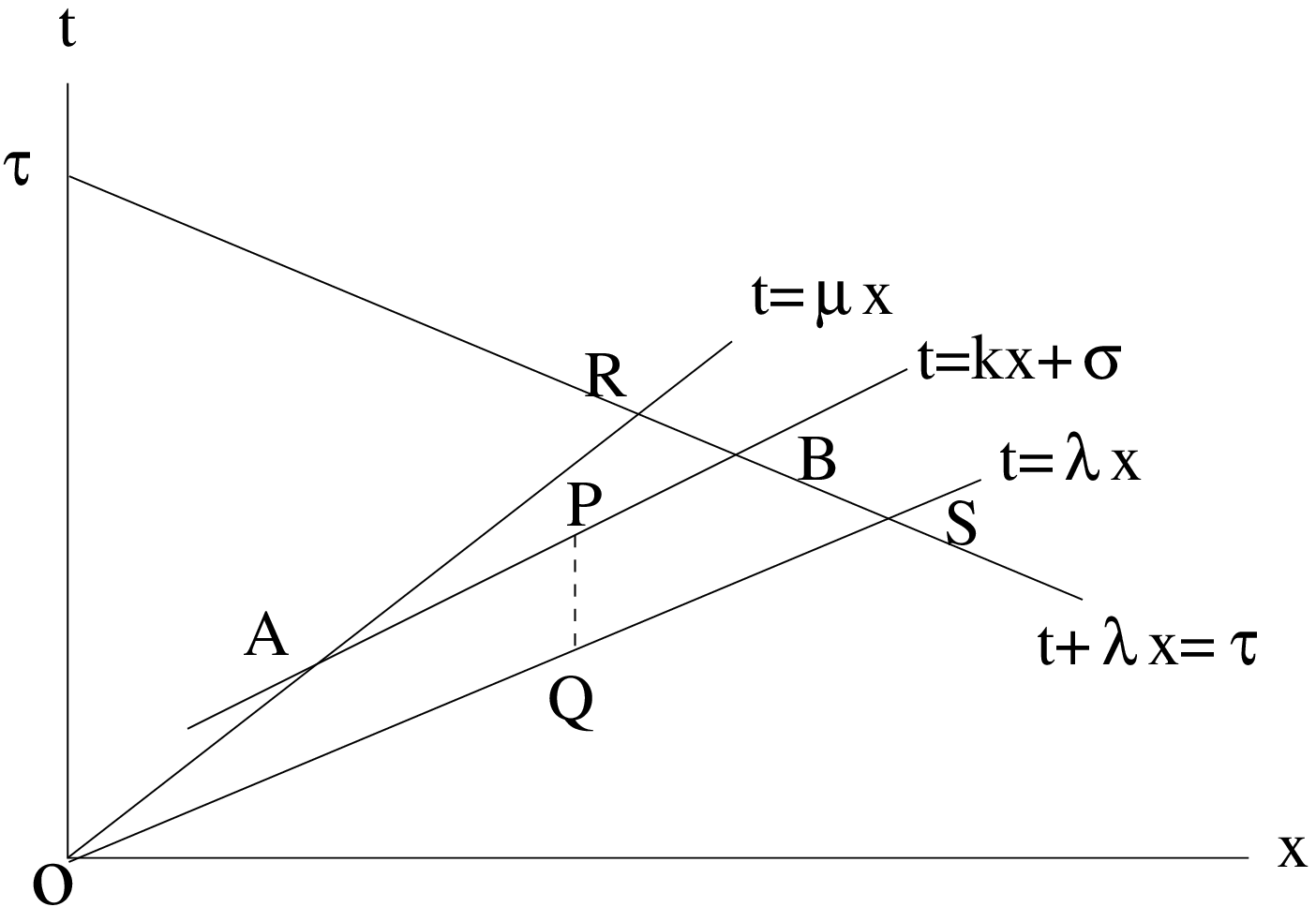, height=1.7in}
\end{center}
\caption{Two speed energy estimates}
\label{fig:slow}
\end{figure}

Fix a slope $k$ strictly between $\lambda$ and $\mu$ and choose an
arbitrary $\sigma \in [\tau - (\lambda+k)x_f , \tau - (\lambda+k)
x_m]$ (the line $t=kx + \sigma$ cuts $t=\tau$ between $R$ and
$S$). In Figure
\ref{fig:slow},
 for certain $\sigma$, $A$ will lie on OS instead of
OR - the calculations below are simpler in this case. We have the
identities
\begin{align}
2 (c^2 w_{tt} - w_{xx})w_t  & = (c^2 w_t^2 + w_x^2)_t - 2 (w_t
w_x)_x,
\label{eq:wt}\\
2 (c^2 w_{tt} - w_{xx})w_x  & = 2(c^2 w_x w_t)_t - (c^2 w_t^2 +
w_x^2)_x.
\label{eq:wx}
\end{align}
Using (\ref{eq:wx})with $w=v_1$, $c=\lambda$ and (\ref{eq:v1char})
and that $v_1$ (and hence $\lambda v_{1t} - v_{1x}$) is zero on
RS, we obtain
\begin{align*}
\iint_{OABS} & 2(\lambda^2 v_{1tt} -v_{1xx})v_{1x} \; dx \; dt
\\
= & \int_{AB} 2 \lambda^2 v_{1x} v_{1t} + k (\lambda^2 v_{1t}^2 +
v_{1x}^2) \; dx - \int_{OS} 2 \lambda^2 v_{1x} v_{1t} +
\lambda (\lambda^2 v_{1t}^2 + v_{1x}^2) \; dx
\\
& + \int_{OA} 2 \lambda^2 v_{1x} v_{1t} +
\mu (\lambda^2 v_{1t}^2 + v_{1x}^2) \; dx
- \int_{BS} \lambda (\lambda v_{1t} - v_{1x})^2 \; dx
\\
= &\int_{AB} \lambda (\lambda v_{1t} + v_{1x})^2 + (k-\lambda)
(\lambda^2 v_{1t}^2 + v_{1x}^2)
\; dx
+ \int_{OA} \lambda (\lambda v_{1t} + v_{1x})^2 + (\mu-\lambda)
(\lambda^2 v_{1t}^2 + v_{1x}^2) \; dx.
\end{align*}
Hence
\beqn
(k-\lambda) \int_{AB}(\lambda^2 v_{1t}^2 + v_{1x}^2)\; dx
\leq
\iint_{OABS}  2(\lambda^2 v_{1tt} -v_{1xx})v_{1x} \; dx \; dt.
\label{eq:v1deriden}
\eeqn
Next, use (\ref{eq:wt}) with $w=v_2$, $c=\mu$ and
(\ref{eq:v1char}); also construct positive $a,b$ with $ab=k$ and
$a<\mu$ and $b<1$. Then we have
\begin{align*}
\iint_{OABS} & 2(\mu^2 v_{2tt} -v_{2xx})v_{2t} \; dx \; dt
\\
= & \int_{AB} \mu^2 v_{2t}^2 + v_{2x}^2 + 2k v_{2t} v_{2x} \; dx -
\int_{OS} \mu^2 v_{2t}^2 + v_{2x}^2 + 2 \lambda v_{2x} v_{2t} \; dx
 \\
 & + \int_{OA} \mu^2 v_{2t}^2 + v_{2x}^2 + 2 \mu v_{2x} v_{2t} \; dx
 + \int_{BS} \mu^2 v_{2t}^2 + v_{2x}^2 - 2 \lambda v_{2t} v_{2x} \; dx
 \\
 = & \int_{AB} (a v_{2t} + bv_{2x})^2 + (\mu^2 - a^2) v_{2t}^2 +
 (1-b^2) v_{2x}^2 \; dx
 \\
 & + \int_{OA} (\mu v_{2t} + v_{2x})^2 \; dx +
 \int_{BS} (\mu^2 - \lambda^2) v_{2t}^2 +  (\lambda v_{2t} - v_{2x})^2 \; dx.
\end{align*}
Hence
\beqn
\int_{AB} (\mu^2 - a^2) v_{2t}^2 + (1-b^2) v_{2x}^2 \; dx
\leq
\iint_{OABS}  2(\mu^2 v_{2tt} -v_{2xx})v_{2t} \; dx \; dt.
\label{eq:v2deriden}
\eeqn
Finally, for $i=1,2$ and any $P$ on $AB$, using (\ref{eq:v1char})
we have
\[
v_i^2(P)  = v_i^2(Q) + 2 \int_{QP} v_i v_{it} \; dt
\leq \int_{QP} v_i^2 + v_{it}^2 \; dt.
\]
Hence
\beqn
\int_{AB} v_i^2 \; dx \leq \iint_{OABS} v_i^2 + v_{it}^2 \; dx \;
dt.
\label{eq:venergy}
\eeqn
If we define
\[
E(\sigma) := \int_{AB} |\vb_t|^2 + |\vb_x|^2 + |\vb|^2 \; dx
\]
then (\ref{eq:v1deriden}), (\ref{eq:v2deriden}),
(\ref{eq:venergy}) may be combined to show that
\begin{align*}
E(\sigma) & \cleq \iint_{OABS} |\vb|^2 + |\vb_x|^2 + |(a_{1}
v_{2x} + b_{11}v_1 +b_{12}v_2)v_{1x}| + |(a_{21} v_{1x} + b_{21}
v_1+b_{22}v_2)v_{2t}|
\; dx
\; dt
\\
& \cleq \iint_{OABS} |\vb|^2 + |\vb_x|^2 + |\vb_t|^2 \; dx \; dt
\\
& \cleq \int_{\tau - (\lambda+k)x_f(\tau)}^\sigma E(s) \; ds.
\end{align*}
Hence $E(\sigma)=0$ on $[\tau - (\lambda+k)x_f(\tau) , \tau -
(\lambda+k) x_m(\tau)]$ by Gronwall's inequality.

\hfill{\bf QED}


\section{Proof of Theorem \ref{thm:progwave}}\label{sec:progwave}
We first prove the uniqueness. If there are two solutions of
(\ref{eq:Upde}), (\ref{eq:Uic}), (\ref{eq:uubBC}) in
$C^2([0,\infty), \cD'(\R))$ then their difference is also a
solution in $C^2([0,\infty), \cD'(\R))$ of (\ref{eq:Upde}),
(\ref{eq:Uic}) but with $\ub(0,t)=0$, $\ubb(0,t)$. Convolving this
difference with any compactly supported smooth function of $t$, we
have a smooth solution of this homogeneous initial boundary value
problem and hence it will be zero by standard energy estimates.
Since the convolution was with an arbitrary function of $t$, the
difference of the two solutions must be zero proving the
uniqueness part of Theorem \ref{thm:progwave}.

If we can construct $\ub$ and $\ubb$ in the forms (\ref{eq:ub})
and (\ref{eq:ubb}) with $\fb, \gb$ being of $C^2$ regularity then
$\ub$ and $\ubb$ will be in $C^2([0, \infty), \cD'(\R))$. Now we
construct expansions for $\ub$ and $\ubb$ with the properties
mentioned in Theorem \ref{thm:progwave}.

\subsection{Progressing wave expansion}\label{subsec:exp}

If $c$ is a constant, $\hb(x,t)$ an arbitrary function and
$s(\cdot)$ a distribution then one may show that
\beqn
\calL \{ \hb s(\tcx ) \} = s(\tcx) \calL \hb
 + s'(\tcx) \calT(\calL, c)\hb
+ \partial_t(s'(\tcx ) \left \{ C - c^2 I \right \} \hb)
 \label{eq:LPf} ~.
 \eeqn
where the first order transport operator $\calT$  is defined as
\beqn
 \calT(\calL, c) \hb :=
 \left \{ (C + c^2 I) \partial_t + 2 c I \partial_x + c A \right \}\hb ~.
 \label{eq:transportdef}
 \eeqn

We seek $\ub$ and $\ubb$ in the form given by (\ref{eq:ub}),
(\ref{eq:ubb}) for some arbitrary $C^2$ functions $\fb, \gb, \fbb,
\gbb$ which we assume are defined for all $x,t$. Of course the
value of $\fb, \gb, \fbb, \gbb$ only on the relevant parts will be
needed to determine $\ub, \ubb$.

From the boundary condition (\ref{eq:uubBC}) we see that
\beqn
\gb(0,t)=0, ~~ \gbb(0,t)=0.
\label{eq:gbBC}
\eeqn
We now determine the conditions determining $\fb$ and $\gb$. Using
(\ref{eq:ub}) and (\ref{eq:LPf}) we have
\begin{align*}
0 = \calL \ub = & ~- \delta(\tlx) B \eb_1 + \lambda \delta'(\tlx) A \eb_1 \\
    &+ H(\tlx) \calL \fb + \delta(\tlx) \calT(\calL,\lambda) \fb +
\partial_t ( \delta(\tlx) (C - \lambda^2 I)\fb ) \\
   &+  H(\tmx) \calL (\gb - \fb) + \delta(\tmx) \calT(\calL,\mu) (\gb - \fb) +
 \partial_t( \delta(\tmx) (C - \mu^2 I)( \gb -\fb) ).
\end{align*}
So\footnote{One may show that the conditions imposed below are not
just sufficient but also necessary to have $\calL \ub =0$ and
$\calL \ubb=0$} we insist that $\fb$ and $\gb$ satisfy
(\ref{eq:fgpde}); further, on $t=\lambda x$ we insist that
\begin{align}
(C- \lambda^2 I)\fb & = -\lambda A\eb_1 \label{eq:val1} \\
\calT(\calL,\lambda) \fb &= B\eb_1  \label{eq:trans1}
\end{align}
and on $t=\mu x$ we insist that
\begin{align}
(C- \mu^2I)(\gb - \fb) &=0 \label{eq:val2} \\
\calT(\calL,\mu) (\gb - \fb) &= 0. \label{eq:trans2}
\end{align}
Now (\ref{eq:val1}), (\ref{eq:val2}) give
\beqn
f_2(x, \lambda x) = \frac{ \lambda}{\lambda^2 -\mu^2} a_{21}(x),
\qquad f_1(x, \mu x) = g_1(x, \mu x);
\label{eq:f2f1g1}
\eeqn
further (\ref{eq:trans1}) implies
\begin{align}
(2 \lambda^2 f_{1t} + 2 \lambda f_{1x} + \lambda a_{12} f_2)(x,
\lambda x) &= b_{11}(x)
\label{eq:trans3} \\
((\lambda^2 +\mu^2) f_{2t} + 2 \lambda f_{2x}  + \lambda
a_{21}f_1)(x,\lambda x) &= b_{21}(x);
\label{eq:trans4}
\end{align}
and (\ref{eq:trans2}) implies
\begin{align}
((\lambda^2+ \mu^2) (g_{1} - f_{1})_t + 2 \mu(g_{1} - f_{1})_x +
\mu a_{12} (g_2 - f_2))(x, \mu x) &= 0
\label{eq:trans5} \\
(2 \mu^2( g_2 - f_{2})_t + 2 \mu (g_2 - f_{2})_x  +\mu a_{21}(g_1
- f_1))(x,\mu x) & = 0.
\label{eq:trans6}
\end{align}
From (\ref{eq:f2f1g1}) we have  $(g_1-f_1)(x, \mu x)=0$, hence the
transport equation (\ref{eq:trans6}) implies that $(g_2-f_2)(x,
\mu x)$ is constant. But $\gb(0,t)=0$ and $f_2(0,0) = \lambda
a_{21}(0)/(\lambda^2 -\mu^2)$. Hence
\[
(g_2-f_2)(x, \mu x) =  \frac{ \lambda a_{21}(0)}{\mu^2 -\lambda^2}
.
\]
Since $(g_1-f_1)(x, \mu x)=0$, taking its derivative and using it
in (\ref{eq:trans5}), we may conclude that
\[
(g_1 - f_1)_t(x, \mu x) = \frac{ \lambda \mu}{ (\mu^2 -
\lambda^2)^2} a_{21}(0) a_{12}(x).
\label{eq:g1f1mu}
\]
Finally, using (\ref{eq:f2f1g1}) in (\ref{eq:trans3}) we obtain
\[
2 \lambda \frac{d}{dx} ( f_1(x, \lambda x) ) = b_{11}(x) +
\frac{\lambda^2}{\mu^2 - \lambda^2} a_{12}(x) \, a_{21}(x).
\]
Integrating this and using $f_1(0,0) = g_1(0,0)=0$, we obtain
\[
f_1(x, \lambda x) = \frac{1}{2 \lambda} \int_0^x b_{11}(z) \, dz +
\frac{\lambda}{2(\mu^2-\lambda^2)}
\int_0^x a_{12}(z) \, a_{21}(z) \, dz.
\]

We now determine the conditions characterizing $\fbb$ and $\gbb$.
Using (\ref{eq:ubb}) and (\ref{eq:LPf}) we have
\begin{align*}
0 = \calL \ubb = & ~- \delta(\tmx) B \eb_2 + \mu \delta'(\tmx) A \eb_2 \\
    &+ H(\tlx) \calL \fbb + \delta(\tlx) \calT(\calL,\lambda) \fbb +
\partial_t( \delta(\tlx) (C - \lambda^2 I)\fbb) \\
   &+  H(\tmx) \calL (\gbb - \fbb) + \delta(\tmx) \calT(\calL,\mu) (\gbb - \fbb) +
 \partial_t( \delta(\tmx) (C - \mu^2 I)( \gbb -\fbb) ).
\end{align*}
So we insist that $\fbb$ and $\gbb$ satisfy (\ref{eq:fgpde});
further on $t=\lambda x$ we insist that
\begin{align}
(C- \lambda^2 I)\fbb & = 0 \label{eq:val21} \\
\calT(\calL,\lambda) \fbb &= 0  \label{eq:trans21}
\end{align}
and on $t=\mu x$ we insist that
\begin{align}
(C- \mu^2I)(\gbb - \fbb) &=-\mu A \eb_2  \label{eq:val22} \\
\calT(\calL,\mu) (\gbb - \fbb) &= B \eb_2 \label{eq:trans22}
\end{align}
From (\ref{eq:val21}) and (\ref{eq:val22}), we obtain
\beqn
\bar{f}_2(x,\lambda x) =0, ~~ (\bar{g}_1-\bar{f}_1)(x, \mu x) = \frac{\mu}{\mu^2 - \lambda^2} a_{12}(x),
\label{eq:charval}
\eeqn
from (\ref{eq:trans21}) we obtain
\begin{align}
(2 \lambda^2 \bar{f}_{1t} + 2 \lambda \bar{f}_{1x} + \lambda
a_{12} \bar{f}_2)(x, \lambda x) &= 0
\label{eq:trans23}
\\
( (\lambda^2 + \mu^2) \bar{f}_{2t} + 2 \lambda \bar{f}_{2x} +
\lambda a_{21} \bar{f}_1)(x, \lambda x) &= 0,
\label{eq:trans24}
\end{align}
and from (\ref{eq:trans22}) we obtain
\begin{align}
( (\lambda^2 + \mu^2) (\barg_{1} - \barf_{1})_t + 2 \mu(\barg_{1}
- \barf_{1})_x + \mu a_{12} (\barg_2 - \barf_2))(x, \mu x) &=
b_{12}(x)
\label{eq:trans25} \\
(2 \mu^2( \barg_2 - \barf_{2})_t + 2 \mu (\barg_2 - \barf_{2})_x
+\mu a_{21} (\barg_1 - \barf_1))(x,\mu x) &= b_{22}(x).
\label{eq:trans26}
\end{align}
Using (\ref{eq:charval}) in (\ref{eq:trans23}) we conclude that
$\barf_1(x, \lambda x)$ is constant. Now $\gbb(0,0)=0$, so from
(\ref{eq:charval}), $\barf_1(0,0)= a_{12}(0) \mu/(\lambda^2 -
\mu^2)$. Hence
\[
\barf_1(x, \lambda x) = \frac{ \mu}{\lambda^2 - \mu^2} a_{12}(0).
\]
Next using (\ref{eq:charval}) in (\ref{eq:trans24}), we conclude
that
\[
\barf_{2t}(x, \lambda x) = - \frac{ \lambda }{ (\mu^2 - \lambda^2)} a_{21}(x) \barf_1(x, \lambda x)
=
\frac{\lambda \mu}{(\lambda^2 - \mu^2)^2} a_{12}(0) a_{21}(x).
\]
Also, from (\ref{eq:trans26}) and (\ref{eq:charval}) we conclude
that
\[
2 \mu \frac{d}{dx} ( (\barg_2 - \barf_2)(x, \mu x)) = b_{22}(x) +
\frac{\mu^2}{\lambda^2 -\mu^2} a_{12}(x) \, a_{21}(x)
\]
which (with initial conditions) implies that
\[
(\barg_2 - \barf_2)(x, \mu x) = \frac{\mu}{2 (\lambda^2 - \mu^2)}
\int_0^x a_{12}(z)
\, a_{21}(z) \, dz
+ \frac{1}{2 \mu} \int_0^x b_{22}(z) \, dz.
\]

Hence to prove Theorem \ref{thm:progwave}, we have to show that
the initial boundary value problems (\ref{eq:fgpde})-(\ref{eq:g2})
and (\ref{eq:fgbpde})-(\ref{eq:bf1}) have $C^2$ solutions over
appropriate regions. Since $A \in C^2[0,\infty)$ and $B \in C^1[0,
\infty)$ the right hand sides of all the zeroth order boundary and
characteristic conditions are $C^2$ functions and the right hand
sides of the first order characteristic conditions are at least
$C^1$ functions. Further the compatibility conditions area also
satisfied at $(0,0)$. Hence Theorem \ref{thm:progwave} follows
from Proposition \ref{prop:goursat} in subsection
\ref{subsec:bvp}.


\subsection{The Characteristic Boundary Value Problem}\label{subsec:bvp}

Pick a constant $T>0$ and define the upper and lower regions
\begin{align*}
U_T &:= \{ (x,t) \, : \, 0 \leq \mu x \leq t \leq T \},
\\
L_T &:= \{ (x,t) \, : \, 0 \leq  \lambda x \leq t \leq \mu x, ~ t
\leq T\},
\\
D_T &= U_T \cup L_T.
\end{align*}

\begin{prop}\label{prop:goursat}

Suppose $A,B$ are in $C^1[0,T/\lambda]$, $\bp \in C^2[0,T]$, $\bq
\in C^2[0,T/\mu]$, $\br \in C^2[0,T/\lambda]$, $s_1 \in C^1[0,
T/\mu]$, $s_2 \in C^1[0, T/\lambda]$ and satisfy the compatibility
condition at $(0,0)$, that is $\bp(0)-\br(0)=\bq(0)$. Also suppose
that $\Fb, \Gb$ are $C^1$ on $L_T, U_T$ respectively.  Then the
Goursat problem
\beqn
\calL \fb = \Fb ~\text{in}~ L_T,
\qquad
\calL \gb = \Gb ~ \text{in} ~ U_T,
\label{eq:fFgG}
\eeqn
with the boundary conditions
\begin{align}
&\gb(0,t) = \bp(t) ~\text{for} ~ t \in [0,T], \\
& (\gb -\fb)(x, \mu x) = \bq(x), ~~ (g_1 -f_1)_t(x, \mu x) =
s_1(x), ~\text{for} ~0 \leq x \leq T/\mu,
\label{eq:mux} \\
& \fb(x, \lambda x) = \br(x), ~~ f_{2t}(x, \lambda x) = s_2(x),
~\text{for} ~ 0 \leq x \leq T/\lambda,
\label{eq:lax}
\end{align}
has a unique solution with $\fb, \gb$ in $C^2$. Further
\beqn
\|\fb \|_{C^2} + \|\gb \|_{C^2} \cleq \|\bp\|_{C^2} + \|\bq\|_{C^2} +
\|\br\|_{C^2} + \|\bs\|_{C^1} + \| \Fb \|_{C^1} + \| \Gb \|_{C^1}
\eeqn
with the constant determined only by $\|A\|_{C^1}, \|B \|_{C^1}$
and $\lambda, \mu,T$.
\end{prop}

\noindent
{\bf Proof of Proposition \ref{prop:goursat}}\\
The uniqueness follows from the analysis in the $A=0, B=0, \Fb=0,
\Gb=0$ case below. We only give an outline of the proof of the
existence part, highlighting the parts of the proof which are not
standard. First, we explicitly write the solution of
(\ref{eq:fFgG})-(\ref{eq:lax}) for the special case when $A=0,
B=0, \Fb=0, \Gb=0$. Then we use this special solution to reduce
the original problem to the case where $\bp=0, \bq=0, \br=0,
\bs=0$ which we deal with using a Volterra equation approach. For
such problems, existence in $H^1$ may be derived by extending
functions by $0$ and appealing to standard results for IBVP on the
region $x \geq 0, t\geq 0$. However, these results will not give
us the higher regularity of $\fb, \gb$ because the total solution
does not have this higher regularity across $t=\mu x$. Hence one
has to appeal to techniques specialized to the problem under
consideration.

\begin{figure}[h]
\begin{center}
\epsfig{file=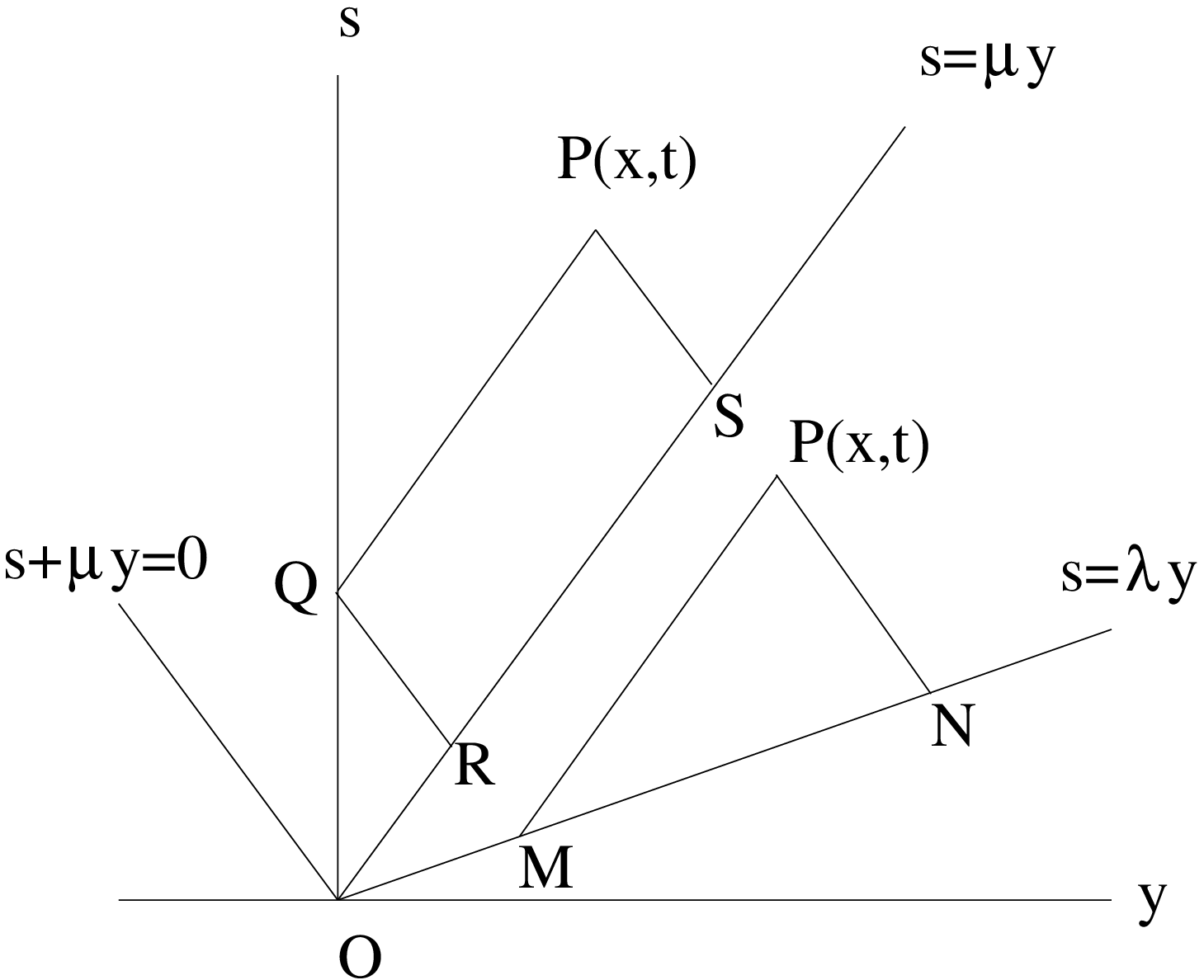, height=1.6in}
\end{center}
\caption{Constructing $f_2, g_2$ when $A=0, B=0, \Fb=0, \Gb=0$ }
\label{fig:homog}
\end{figure}

\noindent
{\bf ($A=0, B=0, \Fb=0, \Gb=0$ case)}\\
In this situation, the equations decouple so the problem reduces
to studying characteristic boundary value problems for the wave
equation. The derivation of the formulas for $f_2, g_2$ is easy
enough; for $(x,t)$ with $\lambda x \leq t \leq \mu x$, $f_2$ is
determined by the values of $f_2$ and $f_{2t}$ on $t=\lambda x$ -
see the triangle $PMN$ in Figure \ref{fig:homog}. Hence, the
transmission condition (\ref{eq:mux}) gives us $g_2$ on $t=\mu x$;
then for any point $P(x,t)$ with $\greg$ we have $g_2(P) = g_2(Q)
+ g_2(S) - g_2(R)$. The expressions for $f_2(x,t)$ and $g_2(x,t)$
consist of the values of $r_2, q_2, p_2$ at linear combinations of
$x,t$ and the integral of $s_2$ over an interval with end points
which are linear combinations of $x,t$. Hence the $C^2$ regularity
of $f_2, g_2$ follows quickly from the regularity of $\bp, \bq,
\br, \bs$.

The derivation of the formula for $f_1, g_1$ is not so clear cut
because now $s= \lambda y$ is a characteristic and hence the value
of $f_1$ on $s=\lambda y$ alone is not enough to determine $f_1$
on $\freg$. An implicit method is needed and the boundary
condition on $x=0$ and the transmission condition on $t=\mu x$ now
play a role. One starts with $f_1$ and $g_1$ as sums of unknown
functions of $\tlx$ and $t+\lambda x$ and the required boundary
and transmission conditions lead to the determination of the
unknown functions. We will not write the long expression for
$f_1(x,t)$ and $g_1(x,t)$ which consist of the values of $r_1,
p_1, q_1$ at linear combinations of $x,t$ and the integral of
$s_1$ over an interval with end points which are linear
combinations of $x,t$. Hence the $C^2$ regularity of $f_1, g_1$
follows from the regularity of $\bp, \bq, \br, \bs$.

\noindent
{\bf (General case)}\\
Let $\pbo(x,t)$ and $\psibo(x,t)$ be the $C^2$ solutions of
(\ref{eq:fFgG})-(\ref{eq:lax}) for the $A=0, B=0, \Fb=0, \Gb=0$
case. Then $\fb - \pbo, \gb - \psibo$ is the solution of
(\ref{eq:fFgG}) - (\ref{eq:lax}) except with $\bp=0, \bq=0, \br=0,
\bs=0$ and $\Fb, \Gb$ replaced by $\Fb+A \pbo_x + B \pbo$, $\Gb +
A \psibo_x + B \psibo$ which are still $C^1$ functions on $L_T$
and $U_T$ respectively. Since $\pbo(x,t)$ and $\psibo(x,t)$ are
$C^2$, we need to prove Proposition \ref{prop:goursat} only for
the case when $\bp=\bq=\br=\bs=0$.

\begin{figure}[h]
\begin{center}
\epsfig{file=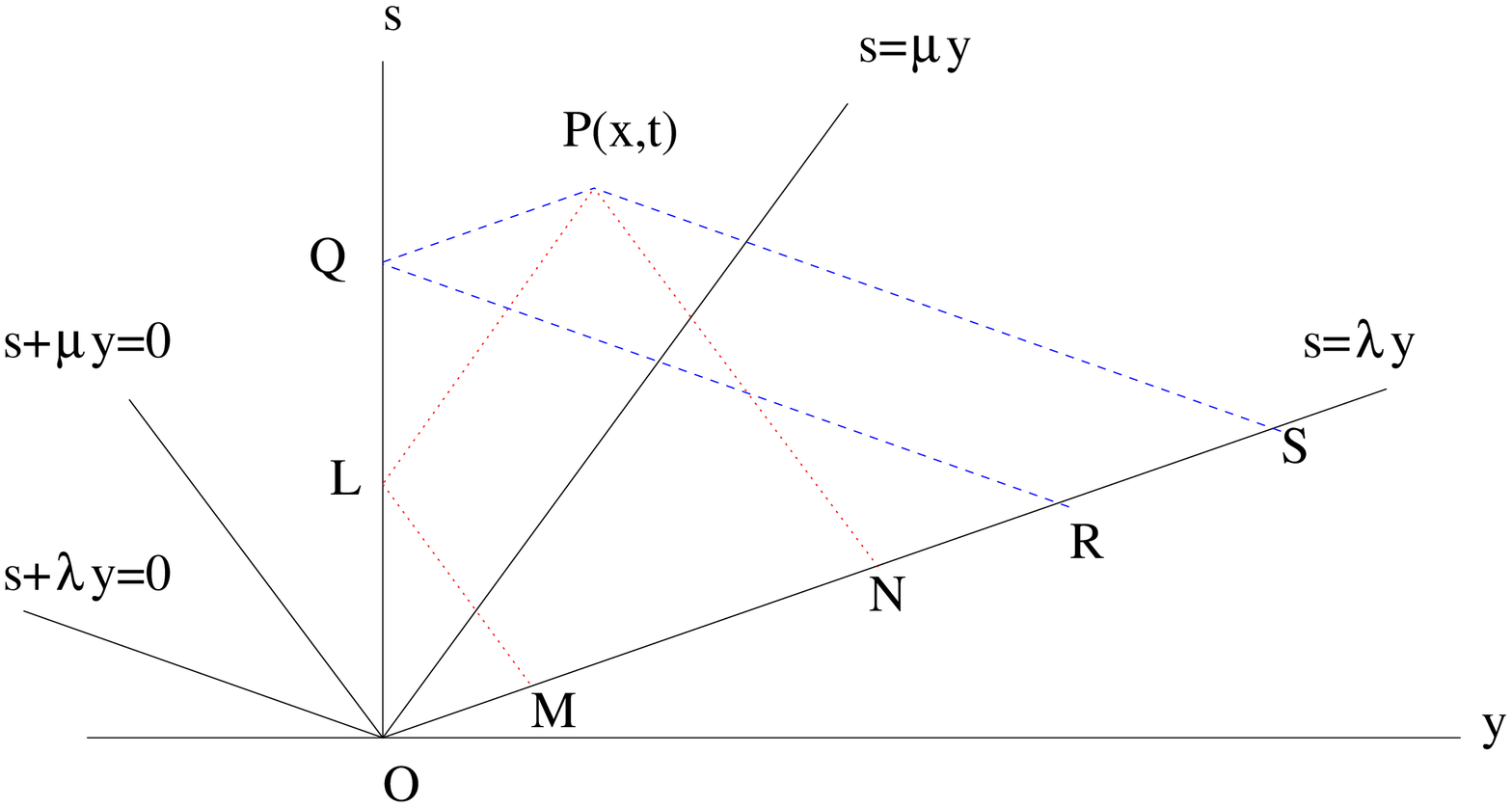, height=1.6in}
\end{center}
\caption{Solution of inhomogeneous wave equation}
\label{fig:inhom}
\end{figure}

For functions $\fb(x,t)$ and $\gb(x,t)$ defined over the regions
$\freg$ and $\greg$ respectively, we define, over the region $0
\leq \lambda x \leq t$, the piecewise function
\[
\{\fb,\gb\}(x,t)
= \begin{cases}
   \fb(x,t) & \text{if} ~\freg, \\
   \gb(x,t) & \text{if} ~ \greg.
   \end{cases}
\]
The value of $\{\fb, \gb\}$ on $t=\mu x$ is ambiguous and is to be
understood to be the one sided limit.
For any vector function $\Hb(x,t)$ on the region $0 \leq \lambda x
\leq t$, we define a vector function $ \I(\Hb)(x,t)$ on the region
$0 \leq \lambda x \leq t$ as (see Figure \ref{fig:inhom})
\[
\I(\Hb)(x,t) = \left [ \iint_{PQRS} H_1, ~ \iint_{PLMN} H_2 \right ]
\]
where PQRS has sides parallel to $s-\lambda y=0$ or $s+\lambda
y=0$ and PLMN has three sides parallel to $s-\mu y=0$ or $s+\mu
y=0$. Note that PLMN will change into a triangle PMN if $ \lambda
x \leq t \leq \mu x$.

Now
\begin{align}
\iint_{PQRS} H_1 & = \frac{1}{2 \lambda} \int_{\tlx}^{t +\lambda x}
\int_0^{\tlx} H_1((q-p)/2\lambda, (q+p)/2) \, dp \, dq,
\label{eq:H1def}
\end{align}
and
\begin{align}
\iint_{PLMN} H_2 &= \frac{1}{2 \mu} \int_{\tmx}^{t+ \mu x}
\int_{ - (\mu-\lambda)q/(\mu+\lambda)}^{\tmx}
H_2((q-p)/2\mu, (q+p)/2) \, dp \, dq
\label{eq:H2def}
\end{align}
with the lower limit of the $q$ integral being
$(\mu+\lambda)|t-\mu x|/(\mu -\lambda)$ in the triangular PMN
case, that is when $t \leq \mu x$.

If $\Hb$ is continuous on the regions $\freg$ and $\greg$ (but may
have jumps across $t= \mu x$) then $\I(\Hb)(x,t)$ is at least
$C^1$ on each of those regions.
Clearly $\I(\Hb)(x,t)$ is continuous on $0 \leq \lambda x \leq t$;
further, the first component of $\I(\Hb)(x,t)$ is $C^1$ on $0 \leq
\lambda x \leq t$ because its derivatives in directions parallel
to $s=\lambda y$ and $s+\lambda y=0$ are the integrals of $H_1$ on
$PS$ and $PQ$ respectively and these vary continuously with $P$
even across $t= \mu x$.
Also, from Figure \ref{fig:inhom}, $\I(\Hb)$ is zero on $x=0$ and
$t=\lambda x$. Also, the first order derivatives of the second
component of $\I(\Hb)$ on $t=\lambda x$ are zero because the
derivatives of the second component, in $\freg$, in directions
parallel $t=\mu x$ and $t+ \mu x= c$ are integrals along $PN$ and
$PM$ (note $L$ is not present in this case). Finally, if $\Hb$ is
$C^1$ on the regions $\freg$ and $\greg$ then $\I(\Hb)(x,t)$ is
$C^2$ on each of those regions. Hence the parts of $\I(\Hb)$ on
$\freg$ and $\greg$ are the unique solution of (\ref{eq:fFgG}) -
(\ref{eq:lax}) when $A=0, B=0$ and $\bp=\bq=\br=\bs=0$ with
$\{\Fb, \Gb\}$ replaced by $\Hb$.

Hence the $\fb, \gb$ we seek are the solutions of the Volterra
like integral equation
\beqn
\{ \fb, \gb \}(x,t) = \I( \{\Fb, \Gb \} )(x,t) +
\I( \{ A \fb_x + B \fb , A \gb_x + B \gb \} )(x,t).
\label{eq:fixed}
\eeqn
Fix a $\bar{T}>0$; we wish to solve (\ref{eq:fixed}) on
$D_{\bar{T}}$.
For any $T \in (0, \bar{T}]$, let $\calB$ be the Banach space of
piecewise vector functions $\{\fb, \gb\}$ with $\fb$ a $C^1$
function on $L_T$ and $\gb$ a $C^1$ function on $U_T$ with the
$C^1$ norms. Now for any $\Hb \in \calB$ we have shown that
$\I(\Hb)$ is also in $\calB$; further using (\ref{eq:H1def}),
(\ref{eq:H2def}), one may show that for $(x,t) \in D_T$
\[
|\I(\Hb)(x,t)|,~ |\partial_t \I(\Hb)(x,t)|, ~ |\partial_x
\I(\Hb)(x,t)|   \cleq (T+T^2) \|\Hb\|_{\calB}.
\]
Here the derivatives of $\I(\Hb)(x,t)$ and $\Hb(x,t)$ are to be
understood to be one sided at points on $t=\mu x$; the constant is
determined only by $\lambda, \mu$. Hence we can define the map
$\K$ from $\calB$ to $\calB$ with
\[
\K(\{\fb,\gb\}) = \I( \{\Fb, \Gb \} )(x,t) +
\I( \{ A \fb_x + B \fb , A \gb_x + B \gb \} )(x,t).
\]
For arbitrary $\{\fb,\gb\}$ and $\{\fb,\gb\}$ in $\calB$, we may
observe that
\[
\| \K( \{\fb,\gb\})
- \K(\{\fb,\gb\}) \|_{\calB} \leq C(T+T^2) ( |A|_{C(D_{\bar{T}})}
+ |B|_{C(D_{\bar{T}})} )
 \|\{\tilde{\fb} -\fb, \tilde{\gb} -\gb\}\|_{\calB}
\]
with $C$ determined by $\lambda$ and $\mu$, implying $\K$ is a
contraction for $T$ small enough. Hence $\K$ has a fixed point in
$\calB$ and we have proved the existence of the unique solution of
(\ref{eq:fixed}) for $T>0$ small enough.

Now suppose we have solved (\ref{eq:fixed}) for $T=T_1$ for some
$0<T_1 < \bar{T}$; so we have $\fb^*, \gb^*$ on $L_{T_1}$ which
solve (\ref{eq:fixed}). For any $T \in [T_1, \bar{T}]$ we redefine
$\calB$ as before except that we require that $\fb=\fb^*$ on
$L_{T_1}$ and $\gb=\gb^*$ on $U_{T_1}$. Define $\K$ as before;
then because $\fb$, $\gb$ agree with $\fb^*$, $\gb^*$ respectively
on $D_{T_1}$ and satisfy (\ref{eq:fixed}) for $T=T_1$, one may see
that
\[
\| \K( \{\tilde{\fb}, \tilde{\gb} \})
- \K(\{\fb,\gb\}) \|_{\calB} \leq C(|T-T_1|+|T-T_1|^2) (
|A|_{C(D_{\bar{T}})} + |B|_{C(D_{\bar{T}})} )
 \|\{\tilde{\fb} -\fb, \tilde{\gb} -\gb\}\|_{\calB}
\]
because after the subtraction and the cancelation of the
contribution to the integrals over the region $D_{T_1}$, the
remaining integrals are over subregions of $T_1 \leq t \leq T$ and
no line parallel to $t = \pm \lambda x$ or $t= \pm \mu x$ in this
region will have a length exceeding a constant times $T-T_1$; the
constant determined by $\lambda, \mu$. Hence as before, $\K$ is a
contraction for $T-T_1$ small enough. The important point is that
there is a positive lower bound on $T-T_1$ for which $\K$ is a
contraction and this lower bound is dependent only on
$|A|_{C(D_{\bar{T}})}$,  $|B|_{C(D_{\bar{T}})}$, $\lambda$, $\mu$
and $\bar{T}$ and is independent of $T_1$. So repeating this
argument we can construct the solution of (\ref{eq:fixed}) over
$D_{\bar{T}}$.

The solution constructed, $\{\fb, \gb\}$, is $C^1$ on $L_T$ and
$U_T$. However for such $\{\fb, \gb \}$ the right hand side of
(\ref{eq:fixed}) is $C^2$ on $L_T, U_T$. Hence $\{\fb, \gb\}$ is
$C^2$ on $L_T$ and $U_T$.

\hfill{\bf QED}


\end{document}